\documentclass{article}
\usepackage[margin=1in]{geometry}
\usepackage[normalem]{ulem}
\usepackage{soul}
\usepackage{amsmath, amsfonts, amssymb, mathtools, enumerate, graphicx, float, color,lscape,tabularx,booktabs}
\newcolumntype{Y}{>{\centering\arraybackslash}X}
\usepackage{algorithm,algpseudocode}
\usepackage[font=footnotesize]{caption}
\usepackage{multirow}

\DeclareMathOperator*{\argmax}{arg\,max}

\title{A Markov decision process approach to optimizing cancer therapy using multiple modalities}
\author{Kelsey Maass$^*$, Minsun Kim$^\dag$ \\ \small{$^*$ Department of Applied Mathematics, University of Washington, Seattle WA, USA} \\ \small{$^\dag$ Department of Radiation Oncology, University of Washington, Seattle WA, USA}}
\date{}

\begin{document}
\maketitle

\abstract
There are several different modalities, e.g., surgery, chemotherapy, and radiotherapy, that are currently used to treat cancer. It is common practice to use a combination of these modalities to maximize clinical outcomes, which are often measured by a balance between maximizing tumor damage and minimizing normal tissue side effects due to treatment. However, multi-modality treatment policies are mostly empirical in current practice, and are therefore subject to individual clinicians' experiences and intuition. We present a novel formulation of optimal multi-modality cancer management using a finite-horizon Markov decision process approach. Specifically, at each decision epoch, the clinician chooses an optimal treatment modality based on the patient's observed state, which we define as a combination of tumor progression and normal tissue side effect. Treatment modalities are categorized as (1) Type 1, which has a high risk and high reward, but is restricted in the frequency of administration during a treatment course, (2) Type 2, which has a lower risk and lower reward than Type 1, but may be repeated without restriction, and (3) Type 3, no treatment (surveillance), which has the possibility of reducing normal tissue side effect at the risk of worsening tumor progression. Numerical simulations using various intuitive, concave reward functions show the structural insights of optimal policies and demonstrate the potential applications of using a rigorous approach to optimizing multi-modality cancer management.

\section{Introduction}

Cancer accounts for nearly 1 in every 7 deaths worldwide, claiming more lives than HIV/AIDS, tuberculosis, and malaria combined \cite{siegel2016cancerfacts}. It has been projected that by 2030 the global burden of cancer will grow to 21.7 million new cancer cases and 13 million cancer deaths worldwide due to growing and aging populations, with more increases expected due to the adoption of behaviors and lifestyles associated with economic development and urbanization \cite{siegel2016cancerfacts}, a trend already observed in economically transitioning countries \cite{torre2015global}. In the United states, cancer is the second leading cause of death, accounting for nearly 1 in every 4 deaths. It is the leading cause of death in 21 states and among adults aged 40 to 79 \cite{siegel2016cancerstats}.  

There are several different modalities that are currently used to treat cancer, including surgery, chemotherapy, radiotherapy, hormone therapy, immune therapy, and targeted therapy \cite{siegel2016cancerfacts, torre2015global}. Most patients receive treatment using two or more modalities, often sequentially, in the course of managing their cancer. As noted in \cite{kemp2016combo}, ``{\it cancer treatments using a single therapeutic agent often result in limited clinical outcomes due to tumor heterogeneity and drug resistance. Combination therapies using multiple modalities can synergistically elevate anti-cancer activity while lowering doses of each agent, hence, reducing side effects.}" Some examples include combining gene therapy and chemotherapy with nanotechnology \cite{kemp2016combo}; using surgery, chemotherapy, and radiotherapy for head and neck cancers \cite{franco2016combined}, combining immunotherapy with the more traditional surgery, radiotherapy, chemotherapy, and targeted therapy \cite{sathyanarayanan2015cancer}; treating glioblastomas with a combination of surgery, radiotherapy, and systemic therapy \cite{taunk2016external}; treating brain metastases using a combination of surgery, radiotherapy, and symptomatic care \cite{matzenauer2016treatment}; treating localized rectal cancer using a combination of surgery, chemotherapy, radiotherapy, and adjuvant cytotoxic therapy \cite{artac2016update}; and chemoradiotherapy plus surgery for esophageal cancer \cite{shapiro2015neoadjuvant}.  

Unfortunately, multi-modality treatment decisions in current practice rely predominantly on individual clinician's experiences, and therefore the optimality of such decisions is unclear. Considering that there are many treatment options available to us in modern medicine and numerous possible outcomes associated with each treatment course, counting on intuition or a heuristic search for optimal multi-modality treatment policies can be costly and inefficient \cite{schaefer2005modeling, bennett2013artificial}. There have been limited efforts to model optimal multi-modality treatments using a mathematical approach. For example, in \cite{beil2001analysis} Beil and Wein studied the optimal sequencing of surgery, chemotherapy, and radiotherapy, using ordinary differential equations to describe the behavior of the primary tumor and metastases. They suggested two novel treatment sequences from their study, but the weakness of their model is the requirement of an accurate, {\it a priori} knowledge of the 14 parameters used in the equations. Alternatively, Hathout et al. investigated the optimal combination of radiotherapy and surgery in the treatment of glioblastoma patients using a reaction-diffusion partial differential equation to simulate the diffusion and proliferation of the tumor and radiation cell-kills \cite{hathout2016modeling}. Their model provides answers to when the extent of surgical resection, in combination with radiotherapy, adds survival benefits.

The purpose of this paper is to propose a novel mathematical framework to optimize multi-modality treatment policies for cancer management using a finite-horizon Markov decision process (MDP) approach, and to demonstrate the feasibility and potential of the proposed model. An MDP with a finite planning horizon is a mathematical framework for optimizing a sequence of actions in a stochastic system, where the state of the system is given at the beginning of each decision epoch. The goal is to maximize the expected reward at the end of the planning horizon. MDPs have been used to model problems in various industries, e.g., robotics \cite{roy1999coastal, kinjo2015evaluation} and economics \cite{briggs1998introduction, rust1996numerical}, and more recently have been successfully applied to problems in medicine. For example, MDPs have been used to find optimal treatments for sperocytosis \cite{magni2000deciding} and ischemic heart disease \cite{hauskrecht2000planning}, to determine the optimal dose in several radiation treatment periods \cite{kim2009markov}, and to decide whether to accept or reject an offered kidney or living-donor liver for transplantation \cite{schaefer2005modeling,ahn1996involving,alagoz2004optimal}. A combination of MDPs and dynamic decision networks has been used to develop a general-purpose artificial intelligence framework that can ``think like a doctor" for personalized medicine, increasing patient outcomes and decreasing costs \cite{bennett2013artificial}. In our model, we define the system state as a combination of the degree of tumor progression (or tumor control) and normal tissue side effect. The action space consists of three different categories of treatment modalities based on the characteristics of repeatability, tumor reduction, and risk to normal tissue.  At each decision epoch the clinician observes the patient's state and chooses an optimal treatment modality accordingly to maximize the expected terminal reward, which is a function of the patient's final state. 

The rest of this paper is organized as follows.  In Section \ref{sec:model} we introduce our model formulation with details on treatment planner's actions (treatment modalities), patient states, state transition probabilities and boundary conditions, and reward functions, where we review the backward induction method to obtain optimal policies. In Section \ref{sec:simulations} we present numerical simulations using various intuitive reward functions to show the structural insights of optimal policies and the potential benefits of using a rigorous model for optimal multi-modality cancer management. Specifically, we demonstrate that changes in the reward functions and state transition probabilities result in changes of the optimal policy that correspond with clinical intuition. Finally, we conclude our paper and discuss possible future extensions to our model in Section \ref{sec:conclusion}. 

\section{Problem formulation}
\label{sec:model}
Consider a treatment course with $T$ periods, where a patient seeks an optimal treatment decision. We define the four components of our MDP model: the treatment planner's actions in Section \ref{subsec:action}, the patient state in Section \ref{subsec:state}, state transition probabilities in Section \ref{subsec:probabilities}, and the intermediate and terminal reward functions in Section \ref{subsec:reward}. Finally, we present the backward induction algorithm to solve the recursive Bellman equations in Section \ref{subsec:problem}. 

\subsection{Treatment planner's actions} 
\label{subsec:action}

We denote the action space as $A = \{M_1, M_2, M_3\},$ where $M_1, M_2,$ and $M_3$ represent the Type 1, 2, and 3 modalities respectively. The definition of each modality type is as follows:

\begin{itemize}
\item Type 1: Treatment modalities with a high risk (increasing side effect) and high reward (decreasing tumor progression). The frequency of administering Type 1 modalities is restricted.
\item Type 2: Treatment modalities with a lower risk and lower reward than Type 1. May be repeated without restriction in frequency.
\item Type 3: No treatment (surveillance). Has a higher probability of reducing normal tissue side effect and increasing tumor progression than the Type 1 and Type 2 modalities. May be repeated without restriction in frequency.
\end{itemize}

We categorize treatment modalities into three types acknowledging that certain modalities may be more effective but limited in the frequency of their administration. For example, whole brain radiotherapy to manage brain metastases is only done once during a patient's lifetime due to the normal tissue side effect associated with it, whereas the gamma knife or partial brain external beam radiotherapy can be administered multiple times until normal tissue toxicity has reached its tolerance level \cite{nccncns}. For such cases, whole brain therapy is categorized as Type 1 while gamma knife and external beam radiotherapy are categorized as Type 2 modalities. 

In general, $M_1$ and $M_2$ may each represent a set of modalities, but to simplify notation we consider the case where there is only one modality of each type in this section. We consider a case with two modalities of Type 2 in Section \ref{subsec:multipleM2}. Similarly, we restrict the Type 1 modality ($M_1$) to one-time use during the course of treatment to simplify notation, though this is certainly not a requirement of our model.

\subsection{Patient state}
\label{subsec:state}
Let the total number of possible states for normal tissue side effect and tumor progression be $m+1$ and $n+1$ respectively.  Let $s_t = (h_t, \phi_t, \tau_t) \in S$ denote the patient state in treatment period $t$ with $t=1,2,\dots, T$.  Each state variable is defined as follows:
\begin{itemize}
\item $h_t \in H=\{0,1\}$: History of Type 1 modality in period $t$ such that $h_t=1$ if the Type 1 modality has been used in periods $1,2,\dots,$ or $t-1$.  Otherwise $h_t=0$.
\item $\phi_t \in \Phi =\{0,1,\dots,m \}$: Normal tissue side effect due to treatment observed in period $t$. $\Phi$ is ordered such that $\phi_t=0$ represents no side effect while $\phi_t=m$ represents the worst possible side effect (patient death due to normal tissue side effect).
\item $\tau_t \in \mathcal{T} =\{0,1,\dots,n \}$: Tumor progression observed in period $t$. $\mathcal{T}$ is ordered such that $\tau_t=0$ represents the best patient state (tumor remission) and $\tau_t=n$ represents the worst possible patient state (patient death due to tumor progression). 
\end{itemize}
We denote the full patient state space by the Cartesian product $S = (H \times \Phi \times \mathcal{T})$.

\subsection{State transition probabilities}
\label{subsec:probabilities}

When treatment modality $a_t \in A$ is implemented in period $t$ for a patient in state $s_t$, the patient's state in period $t+1$ is assumed to be $s_{t+1}$ with the probability of $P_t(s_{t+1}|s_t,a_t)$. We assume that the transition probabilities for each state variable are conditionally independent of one another, that is, 
\begin{equation}
P_t(s_{t+1}|s_t, a_t) = P_t^H(h_{t+1}|h_t,a_t) \times P_t^{\Phi}(\phi_{t+1}|\phi_t,a_t) \times P_t^{\mathcal{T}}(\tau_{t+1}|\tau_t,a_t). 
\end{equation}
This assumption is in line with our intuition that tumor progression and normal tissue side effects do not depend on each other but depend on the treatment type only. The classification of $M_1$ as having a higher risk and higher reward than $M_2$ can then be written in terms of state transition probabilities as follows:
\begin{eqnarray}
 \mbox{Higher risk in side effects:} && P_t^{\Phi}(\phi_{t+1}|\phi_t,M_1) \geq P_t^{\Phi}(\phi_{t+1}|\phi_t,M_2) \mbox{ for } \phi_{t+1} > \phi_t \\
\mbox{Higher reward in tumor control:} &&  P_t^{\mathcal{T}}(\tau_{t+1}|\tau_t,M_1) \geq P_t^{\mathcal{T}}(\tau_{t+1}|\tau_t,M_2) \mbox{ for } \tau_{t+1} < \tau_t
\end{eqnarray}

We define the transition probabilities for the history variable, $h_t$, deterministically to restrict the number of administrations of $M_1$ during the course of treatment as follows:
\begin{eqnarray}
P^H_t(1|h_t,M_1) \hspace{-0.8em}&=&\hspace{-0.8em} 1  \text{ for all } h_t \in H,\\
P^H_t(h_{t}|h_t,a_t) \hspace{-0.8em}&=&\hspace{-0.8em} 1  \text{ for } a_t \in\{M_2,M_3\}. 
\label{eqn:history}
\end{eqnarray}
The restriction that Type 1 modalities may only be used once during the course of treatment is imposed through the transition probabilities in the following manner: 
\begin{equation}
P_t (1,m,n |1,\phi_t,\tau_t,M_1) = 1 \mbox{ for all } \phi_t \in \Phi \text{ and } \tau_t \in \mathcal{T}. 
\end{equation}
This means that the patient will transition to the worst possible state, i.e., $s_{t+1} = (1,m,n)$, if $M_1$ is chosen when the history of $M_1$ use is positive.

We impose absorbing boundary conditions to simulate either the death of the patient (when side effect or tumor progression reaches their maximum value) or tumor remission (when tumor progression reaches zero). We note that in the remission state only tumor progression is fixed, allowing side effect to improve in subsequent treatment periods.
\begin{eqnarray}
\mbox{ Death due to side effect:} && P_t (h_t, m, \tau_t | h_t, m, \tau_t, a_t) = 1 \mbox{ for all } a_t \in A \\
\mbox{ Death due to tumor progression:} && P_t (h_t, \phi_t, n | h_t, \phi_t, n ,a_t ) = 1 \mbox{ for all } a_t \in A \\ 
\mbox{ Tumor remission:} && P_t^{\mathcal{T}}(0  | 0,a_t) = 1 \mbox{ for all } a_t \in A
\end{eqnarray}
	
\subsection{Reward functions}
\label{subsec:reward}
We denote the real-valued terminal reward function $r_{T+1}(s)$, which quantifies the patient's utility of being in state $s$ at the end of their treatment course, $t = T+1$. After each treatment period the patient may also receive an intermediate reward, $r_t(s_t,a_t,s_{t+1})$, which is associated with an action chosen in period $t$, such as the cost of using treatment $a_t$ and the expected outcome in the patient's next state, $s_{t+1}$. While our model will work with any reward function, in general patient utility corresponding to better states should be at least as large as patient utility corresponding to worse states. This means that
\begin{eqnarray}
&r(\phi,\tau) \geq r(\phi',\tau) \mbox{ for } \phi \leq \phi', \\
&r(\phi,\tau) \geq r(\phi,\tau') \mbox{ for } \tau \leq \tau',
\end{eqnarray}
because states are ordered such that smaller states represent better patient conditions.  Some commonly used utility measures in medicine are the quality-adjusted life year, the disability-adjusted life year, and the healthy-years equivalent \cite{weinstein2009qalys, murray1994quantifying, mehrez1989quality,Saokaew2016}. 
	
\subsection{Bellman equations and backward induction}
\label{subsec:problem}

For each patient state $s_t \in S$ and treatment period $t \in \{1,2,\dots,T\}$, our goal is to maximize the expected patient utility at the end of the treatment course, that is, $ E \left[ \sum_{t=1}^{T} r_t (s_t, a_t, s_{t+1}) +r_{T+1} (s_{T+1}) \right]$. Bellman's recursive equations to solve this problem is given by
\begin{equation}
 V_t(s) = \sum_{s' \in S} P_t(s' | s,a) \Big( r_t(s,a,s') + V_{t+1}(s') \Big) \mbox{ for } t=1,2,\c\dots, T
 \end{equation}
with boundary condition $V_{T+1}(s) = r_{T+1}(s)$. The optimal policy can be solved recursively for all $s_t \in S$ and $t = 1,2,\dots,T$ with the well-known backward induction algorithm \cite{puterman2014markov}:

\begin{center}
\hspace{7em}\begin{minipage}{0.7\linewidth}
\begin{algorithmic}
\State{Set $V_{T+1}(s) = r_{T+1}(s)$ for all $s \in S$}
\For{$t = T,T-1,\dots,1$}
	\State{$\displaystyle V_t(s) = \max_{a \in A} \sum_{s'\in S} P_t(s'|s,a) \Big(r_t(s,a,s') + V_{t+1}(s')\Big)$}
	\State{$\displaystyle a_t(s) = \argmax_{a \in A}  \sum_{s'\in S} P_t(s'|s,a) \Big(r_t(s,a,s') + V_{t+1}(s')\Big)$}
\EndFor 
\end{algorithmic}
\end{minipage}
\end{center}

With increasing dimensions in the state space $S$ and the action space $A$, the problem faces the ``curses of dimensionality" and may be solved using approximate dynamic programming \cite{powell2011}.

\section{Numerical simulations}
\label{sec:simulations}

In this section we present numerical simulations to illustrate the structure of optimal multi-modality treatment policies generated with our MDP model. In Section 3.1 we present a base case of the state transition probabilities and reward functions used in our simulations, including the general assumptions made to make them clinically relevant. Next, in Sections \ref{subsec:variousreward} and \ref{subsec:intermediate} we explore how changes in the terminal and intermediate reward functions affect the optimal policies. The effect of changes in transition probabilities on the optimal policies is presented in Section \ref{subsec:varioustransition}. Finally, we explore a scenario with multiple Type 2 treatment modalities in Section \ref{subsec:multipleM2}. We note that the general assumptions made in this section are specific to our numerical simulations and not necessary for solving optimal multi-modality treatment policies using our model. In practice, transition probabilities can be estimated from correlations between treatment modalities and patient outcomes from the clinical literature, which is beyond the scope of this paper. 

\subsection{Base case}
\label{subsec;numstatetransition}

We define a base case with one modality of each type, where $A = \{M_1, M_2, M_3\}$, without intermediate rewards, letting $r_t=0$ for $t=1,2, \dots, T$. In all of our numerical simulations, we use 11 states for $\phi$ and $\tau$, i.e., $m = n = 10$, and three treatment periods, that is, $T=3$. First, in Section \ref{subsub:probabilities} we specify our assumptions used to assign state transition probabilities, and we introduce our reward functions derived from clinical intuition and practice in Section \ref{subsub:reward}. This is followed by the optimal policies resulted from these base transition probabilities and rewards in Section \ref{subsub:baseresults}.

\subsubsection{State transition probabilities}
\label{subsub:probabilities}
For our numerical simulations, we utilize stationary transition probabilities that depend only upon the changes between states rather than on the actual value of the state. This means that
\begin{eqnarray}
P_t(s_{t+1}|s_t,a_t) \hspace{-0.8em}&=&\hspace{-0.8em} P(s_{t+1}|s_t,a_t) \mbox{ for } t = 1,2,\dots, T, \\
P^{\Phi}(\phi_{t+1}|\phi_t,a_t) \hspace{-0.8em}&=&\hspace{-0.8em} P^{\Phi}(\phi_{t+1}'|\phi_t',a_t)  \mbox{ whenever } \phi_{t+1}-\phi_t = \phi_{t+1}' - \phi_t' , \text{ and }\\
P^{\mathcal{T}}(\tau_{t+1}|\tau_t,a_t) \hspace{-0.8em}&=&\hspace{-0.8em} P^{\mathcal{T}}(\tau_{t+1}'|\tau_t',a_t)  \mbox{ whenever } \tau_{t+1}-\tau_t = \tau_{t+1}' - \tau_t'.
\end{eqnarray}

For simplicity, we assume that state variables can only change by one increment between two successive treatment periods, and that tumor progression only improves after treatment ($M_1$ or $ M_2$) while side effect only improves after surveillance ($M_3$). Therefore, tumor progression can either stay the same or get better (decrease by one) between two successive treatment periods after $M_1$ or $M_2$ is chosen, while side effect can either stay the same or get worse (increase by one). This implies that 
\begin{equation}
P^{\mathcal{T}}(0 | 0, a_t) = 1 \mbox{ and } P^{\Phi}(m | m,a_t) = 1 \mbox{ for } a_t \in \{M_1,M_2\},
\end{equation}
which is consistent with our absorbing boundary conditions. When $M_3$ (surveillance) is chosen, tumor progression can either stay the same or get worse between two successive treatment periods, while side effect can either stay the same or get better. This implies that 
\begin{equation}
P^{\mathcal{T}}(n | n,M_3) = 1 \mbox{ and } P^{\Phi}(0 | 0,M_3) = 1. 
\end{equation}

For our numerical simulations, we use the state transition probabilities for non-boundary states defined in Table \ref{tab:baseprobability}. The transition probabilities for the history variable $h$ is as defined in Equation (\ref{eqn:history}). 

\begin{table}[H]
\captionsetup{width=0.9\textwidth}
\centering
\begin{tabularx}{0.9\textwidth}{c *{6}{Y}}
\toprule
\multirow{2}{*}{\bf Modality} 
& \multicolumn{3}{c}{\bf Side Effect in period $t+1$}  
& \multicolumn{3}{c}{\bf Tumor Progression in period $t+1$}\\
\multirow{2}{*}{$(a_t)$} 
& \multicolumn{3}{c}{$(\phi_{t+1})$} 
& \multicolumn{3}{c}{$(\tau_{t+1})$} \\ 
\cmidrule(lr){2-4} \cmidrule(l){5-7}
& $\phi_t-1$ & $\phi_t$ & $\phi_t+1$ & $\tau_t-1$ & $\tau_t$ & $\tau_t+1$ \\
\midrule
$M_1$ &    0 & 0.4 & 0.6 & 0.7 & 0.3 &    0 \\
$M_2$ &    0 & 0.6 & 0.4 & 0.6 & 0.4 &    0 \\
$M_3$ & 0.6 & 0.4 &    0 &    0 & 0.3 & 0.7 \\
\bottomrule
\end{tabularx}
\caption{State transition probabilities $P^{\Phi}(\phi_{t+1}|\phi_t,a_t)$ and $P^{\mathcal{T}}(\tau_{t+1}|\tau_t,a_t)$ used in the base case.}
\label{tab:baseprobability}
\end{table}
	 
\subsubsection{Reward functions}
\label{subsub:reward}

The reward (patient utility) of being in state $s$ depends on the patient's side effect ($\phi$) and tumor progression ($\tau$).  We define our patient utility using additively separable reward functions, where $f$ measures the utility of the normal tissue side effect being in state $\phi$, and $g$ measures the utility of tumor progression being in state $\tau$. Let $d_\phi$ and $d_\tau$ be the parameters associated with functions $f$ and $g$.  Then the total reward is defined as
\begin{equation}
r(\phi,\tau) = c_\phi f(\phi;d_\phi) + c_\tau g(\tau;d_\tau),\\[1ex]
\end{equation}
where $c_\phi$ and $c_\tau$ are the weighting factors of $f$ and $g$, which represent the relative importance of the side effect and tumor progression respectively.  

We use concave reward functions for $f$ and $g$, where improvements made in worse patient states are more appreciated than improvements made in healthier states.  These relationships can often be found in studies that examine patient utility as a function of clinical states \cite{attema2016elicitation,miravitlles2015clinical, carradice2011modelling,currie2006multivariate}. Thus, the functions $f$ and $g$ can be written as
\begin{eqnarray}
\text{Side effect:} \quad f(\phi;d_\phi) \hspace{-0.8em}&=&\hspace{-0.8em} \frac{100}{m^{d_\phi}} \left(m^{d_\phi} - \phi^{d_\phi}\right), \\
\text{Tumor progression:} \quad g(\tau;d_\tau) \hspace{-0.8em}&=&\hspace{-0.8em} \frac{100}{n^{d_\tau}}\left(n^{d_\tau} - \tau^{d_\tau}\right),
\end{eqnarray}
where $d_\phi \geq 1$ and $d_\tau \geq 1$ to make $f$ and $g$ concave.  Note that $d_\phi=d_\tau=1$ indicates a linear function. The  functions $f$ and $g$ are normalized to 100 so that the reward ranges between 0 (minimum) and 100 (maximum).  We note that the concave shape of the reward function is not necessary to our model, but is utilized in our examples to simulate clinically relevant scenarios as previously demonstrated in \cite{kim2009markov}.

Finally, we use a stationary intermediate reward function that does not depend on the action, i.e., $r_t(s_t,a_t,s_{t+1}) = r(s_t, s_{t+1})$ for $t = 1,2,\dots, T$.  For the base case, we use $d_\phi = d_\tau = 2$ and $c_\phi=c_\tau = 1/2$ in the terminal reward function ($r_{T+1}$), without intermediate reward functions, that is, $r (s_t,s_{t+1}) = 0$ for $t=1,2,\dots,T$. 

\subsubsection{Numerical results of the base case}
\label{subsub:baseresults}

\begin{figure}[H]
\captionsetup{width=0.9\textwidth}
\centering
\includegraphics[height=2in]{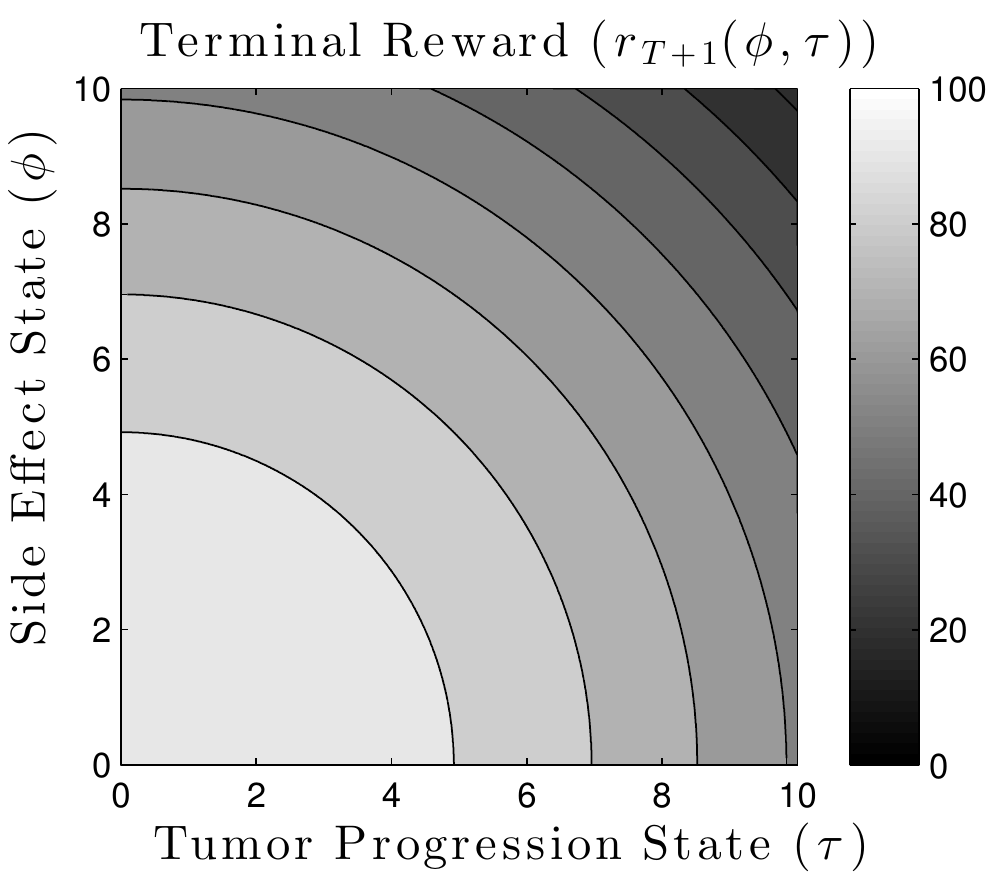} \quad
\includegraphics[height=2in]{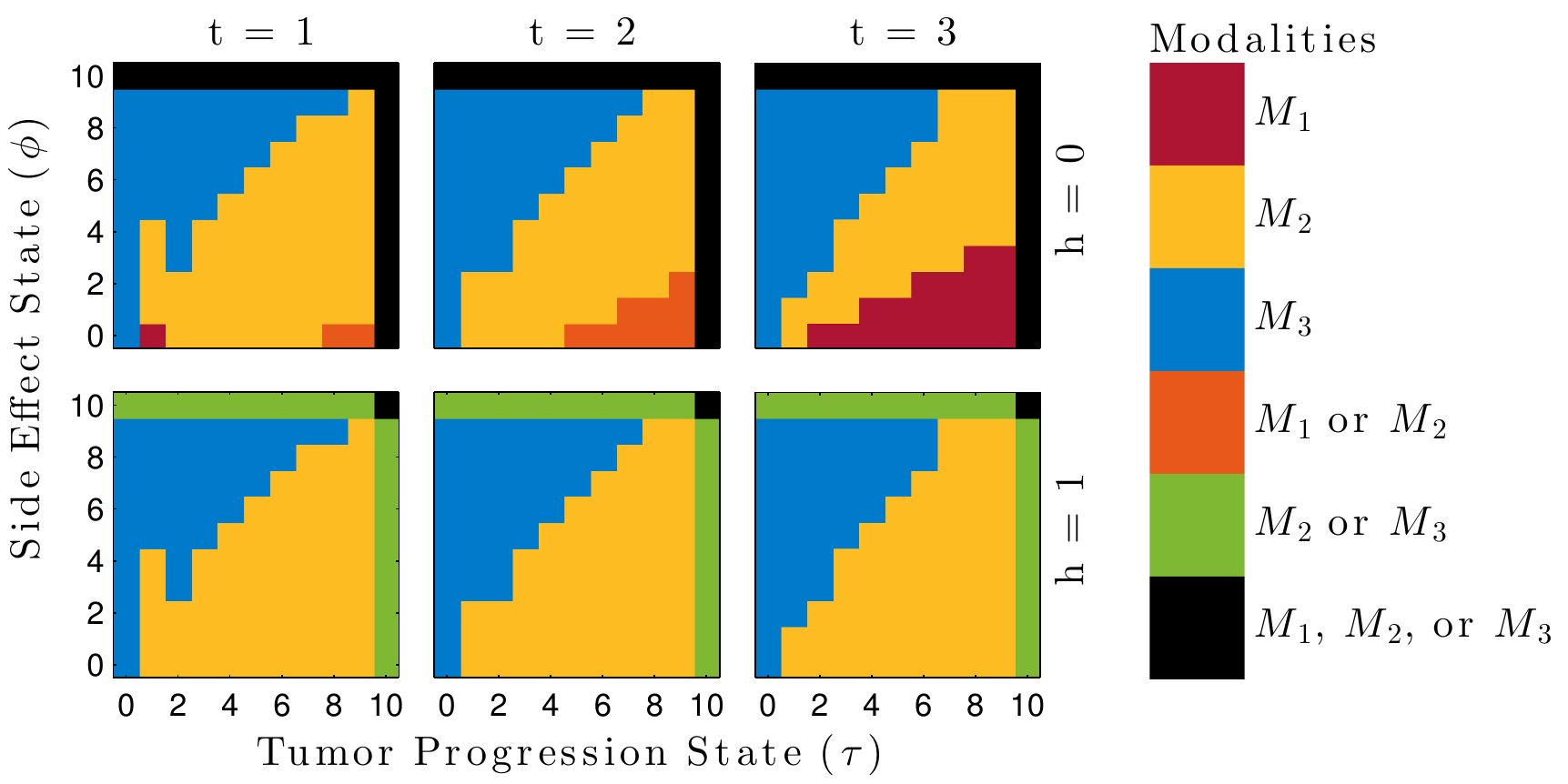}
\caption{(Left) Quadratic terminal reward function, $r_{T+1}(\phi,\tau) = \frac{1}{2} f(\phi; 2) + \frac{1}{2} g(\tau; 2)$. (Right) Optimal treatment policy for $T=3$ and $A=\{M_1, M_2, M_3\}$ with transition probabilities given in Table 1.}
\label{fig:basecase}
\end{figure}

Figure \ref{fig:basecase} shows the optimal policy for the base case with transition probabilities from Table \ref{tab:baseprobability} and rewards functions described in Section \ref{subsub:reward}. The terminal reward function is shown on the left of the figure, where the contours show the isolines of the reward for each given state defined with tumor progression on the horizontal axis and normal tissue side effect on the vertical axis. With lighter regions denoting larger rewards, we see that the lower left corner is the patient's best state, i.e. zero side effect and tumor remission, while the upper right corner is the patient's worst state, i.e. maximum side effect and tumor progression. On the right of the figure, the optimal policy is shown for states defined with tumor progression on the horizontal axis and normal tissue side effect on the vertical axis. The top row shows the optimal policy when the Type 1 modality has never been used, i.e. $h=0$, and the bottom row shows the optimal policy for $h=1$. Each column represents a distinctive treatment period.

The results agree with clinical intuition that, in general, when the patient's normal tissue side effect is more detrimental to their overall health than their tumor progression, i.e., the patient has severe side effects but their tumor progression is less concerning (upper left corner), the optimal action is surveillance ($M_3$). On the other hand, when the patient's current state is on the lower right corner, that is, severe tumor progression with minimal side effects, the optimal action is to use more aggressive treatment modalities ($M_1$ if it has not been used). At the absorption states $\phi = m$ or $\tau = n$, there is no difference among modalities when $M_1$ has not been used. However, since using $M_1$ when $h = 1$ deterministically brings the patient's next state to the worst possible state, the optimal policy with $h = 1$ does not include $M_1$ for any states or treatment periods except for the worst state, $s_t = (1,m,n)$. We note that with $h=0$ the Type 1 modality tends to be saved for later treatment periods, so $M_2$ is used more in the beginning of the treatment course. This observation will be compared with the case in Section \ref{subsec:intermediate} when an intermediate reward is added, where the patient's state during intermediate treatment periods contributes to the total reward, unlike the base case where the total reward depends only on the patient's final state.

\subsection{Effect of terminal reward function shapes on optimal policy}
\label{subsec:variousreward}

In this section we demonstrate how changes in the terminal reward function affect optimal policies. Specifically, we consider various shapes of the terminal reward function in Section \ref{subsubsec:shapes} and the relative importance of side effect and tumor progression in the terminal reward function in Section \ref{subsubsec:importance}.

\subsubsection{Shape of terminal reward functions}
\label{subsubsec:shapes}

We examine how the shape of the terminal reward function, represented by the exponents $d_\phi$ and $d_\tau$, affects the optimal policy when side effect and tumor progression are weighted equally, that is, $c_\phi=c_\tau=1/2$. We use the same exponent for $f$ and $g$, letting $d_\phi=d_\tau \equiv d$.  Therefore, we compare the optimal policies generated using the following terminal reward function:
\begin{eqnarray}
 r_{T+1}(\phi,\tau) \hspace{-0.8em}&=&\hspace{-0.8em} \frac{1}{2} f(\phi,d) + \frac{1}{2} g(\tau,d) \\
 \hspace{-0.8em}&=&\hspace{-0.8em} \frac{50}{10^d} \left\{(10^d - \phi^d) + (10^d - \tau^d) \right\}, \\
\mbox{where } d \hspace{-0.8em}&\in &\hspace{-0.8em} \{3/2,2 \mbox{(base case)}, 3\}.
 \end{eqnarray}

\begin{figure}[H]
\captionsetup{width=0.9\textwidth}
\centering
\includegraphics[height=2in]{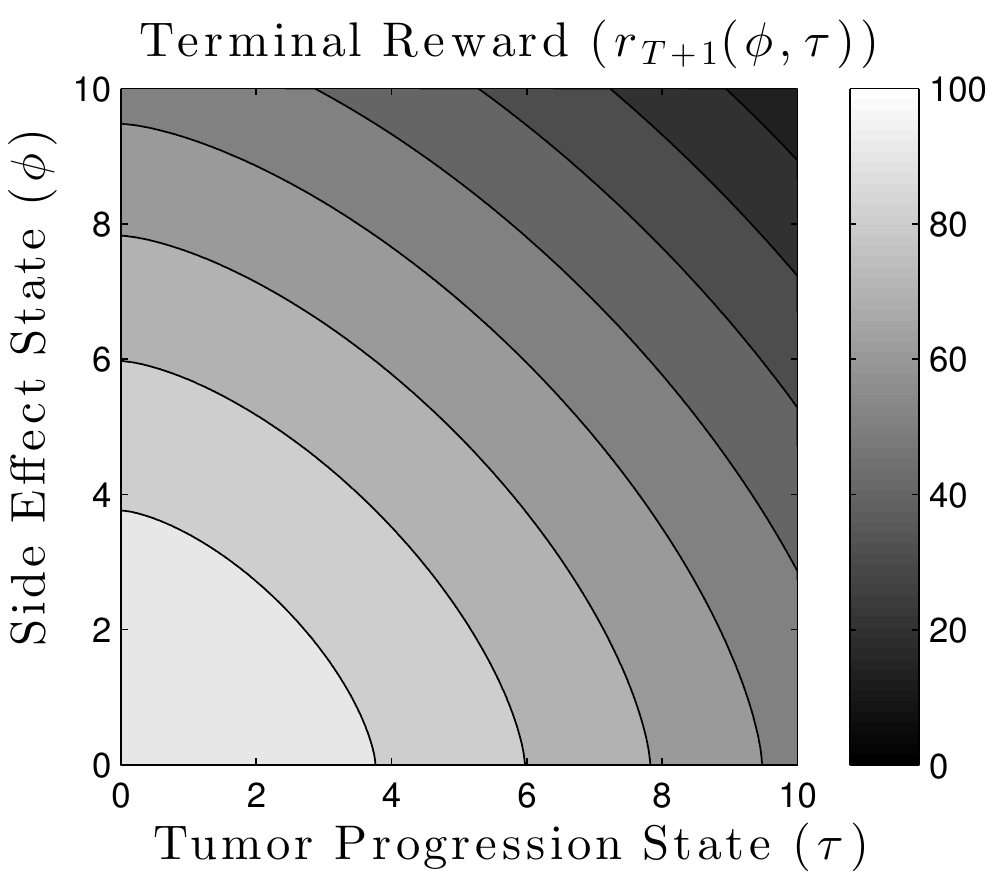} \quad
\includegraphics[height=2in]{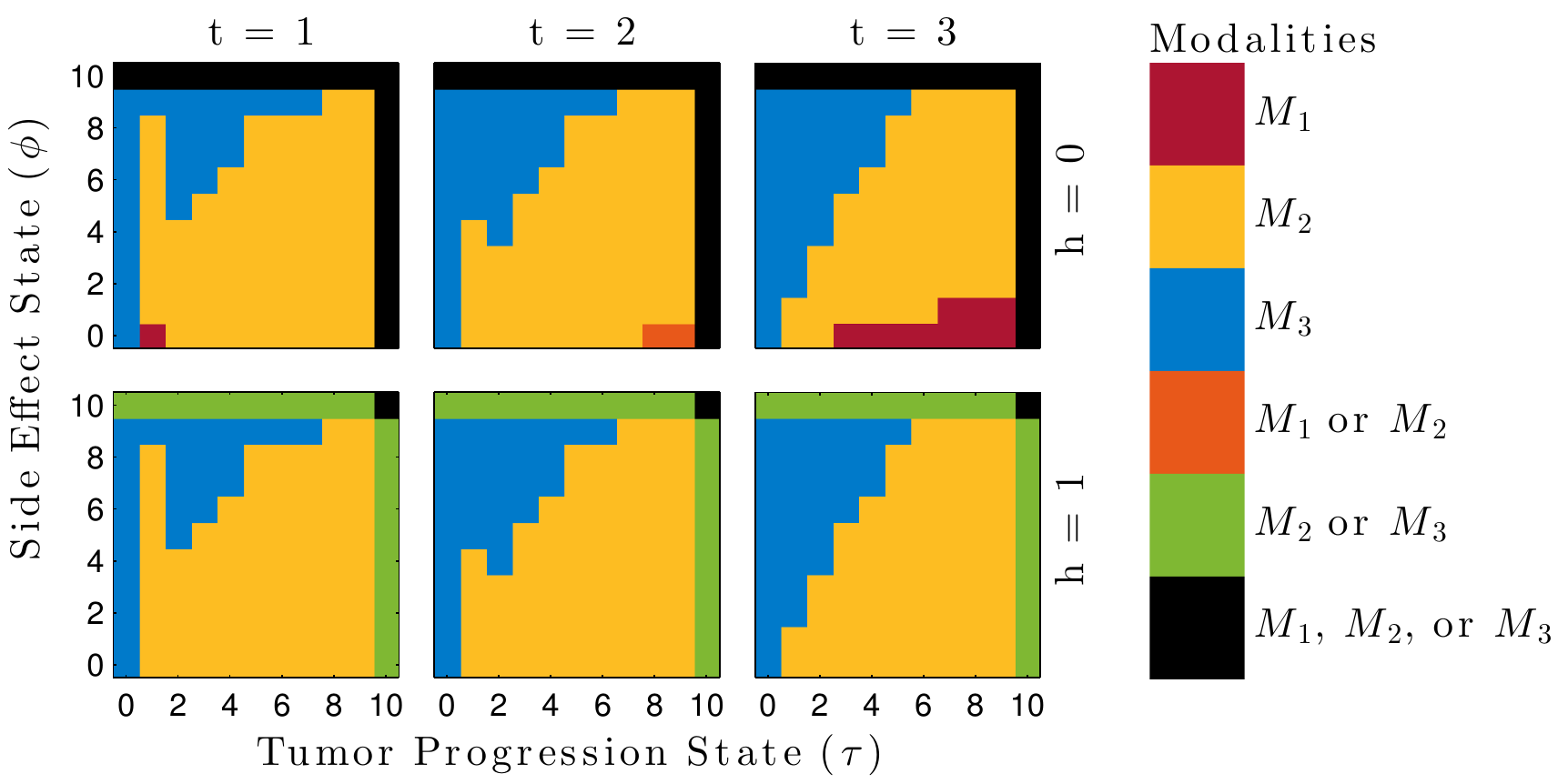}
\caption{(Left) Terminal reward function $r_{T+1}(\phi, \tau) = \frac{1}{2} f(\phi;3/2) + \frac{1}{2} g(\tau;3/2)$. (Right) Optimal treatment policy for $T=3$ and $A=\{M_1, M_2, M_3\}$ with base transition probabilities given in Table \ref{tab:baseprobability}.}
\label{fig:smallerexponent}
\end{figure}

\begin{figure}[H]
\captionsetup{width=0.9\textwidth}
\centering
\includegraphics[height=2in]{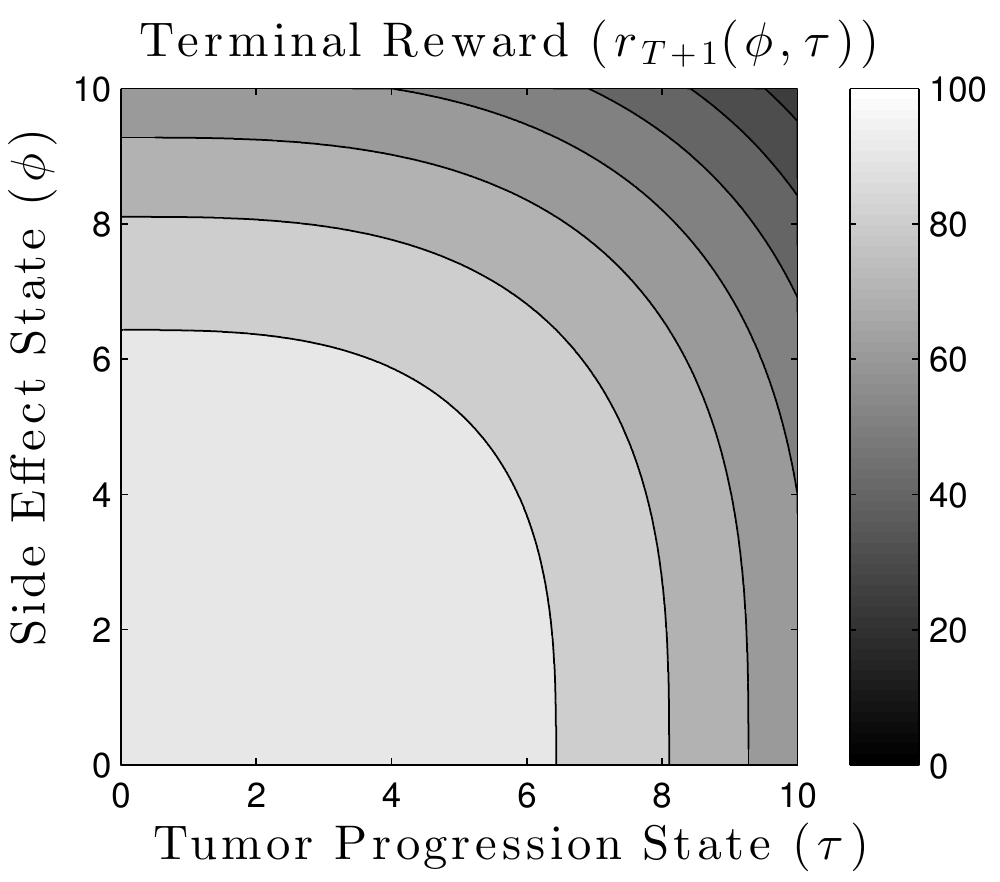} \quad
\includegraphics[height=2in]{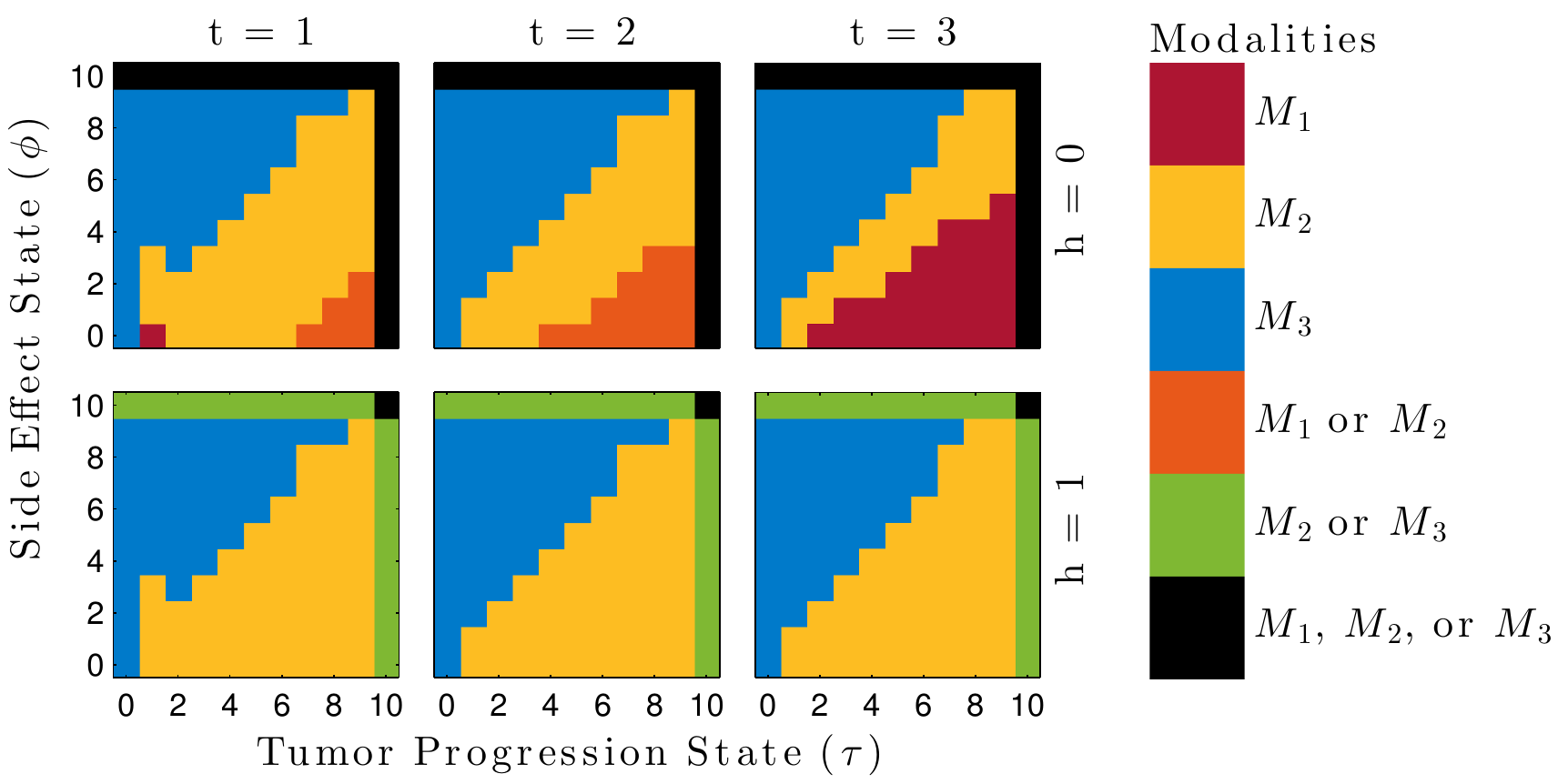}
\caption{(Left) Cubic terminal reward function, $r_{T+1}(\phi,\tau) = \frac{1}{2} f(\phi;3) + \frac{1}{2} g(\tau;3)$.  (Right) Optimal treatment policy for $T=3$ and $A=\{M_1, M_2, M_3\}$ with base transition probabilities given in Table \ref{tab:baseprobability}.}
\label{fig:largerexponent}
\end{figure}

The exponent $d=1$ corresponds to a linear function, where the reward of decreasing tumor progression or side effect is identical in all patient states.  When $d>1$, the reward received from decreasing tumor progression or side effect in worse patient states is larger than the reward received in better states. As $d$ increases, the reward at each patient state also increases, except at the best and worst states where rewards are fixed at zero and 100. This in turn results in a steeper slope and faster transition from ``bad" states to ``good" states. The optimal policies resulting from $d=3/2$ and $d=3$ are shown in Figure \ref{fig:smallerexponent} and Figure \ref{fig:largerexponent} respectively, which can be compared with the base case in Figure \ref{fig:basecase}.
	
As $d$ gets closer to unity, the reward function becomes linear, so improvements made in any patient state (better or worse) are rewarded the same way. Therefore, Figure \ref{fig:smallerexponent} with $d=3/2$ shows that the optimal policy in the upper left corner (high side effect) and the lower right corner (worse tumor progression) is more like the rest of the states as compared to the base case.   As we increase the exponent $d$, the improvements made near the worst patient states (the worst tumor progression or worst side effect) are rewarded higher than the improvements made in other states. This produces treatments that prioritize getting the patient out of these states, with both more surveillance in states with high side effect (upper left corner) and more of the $M_1$ modality in states with worse tumor progression (lower right corner). The optimal policy with $d=3$ is shown in Figure \ref{fig:largerexponent}.  For the remainder of our numerical simulations, we use quadratic terminal reward functions, letting $d = 2$, in order to compare with the base case.

\subsubsection{Relative importance of side effect and tumor progression}
\label{subsubsec:importance}

\begin{figure}[H]
\captionsetup{width=0.9\textwidth}
\centering
\includegraphics[height=2in]{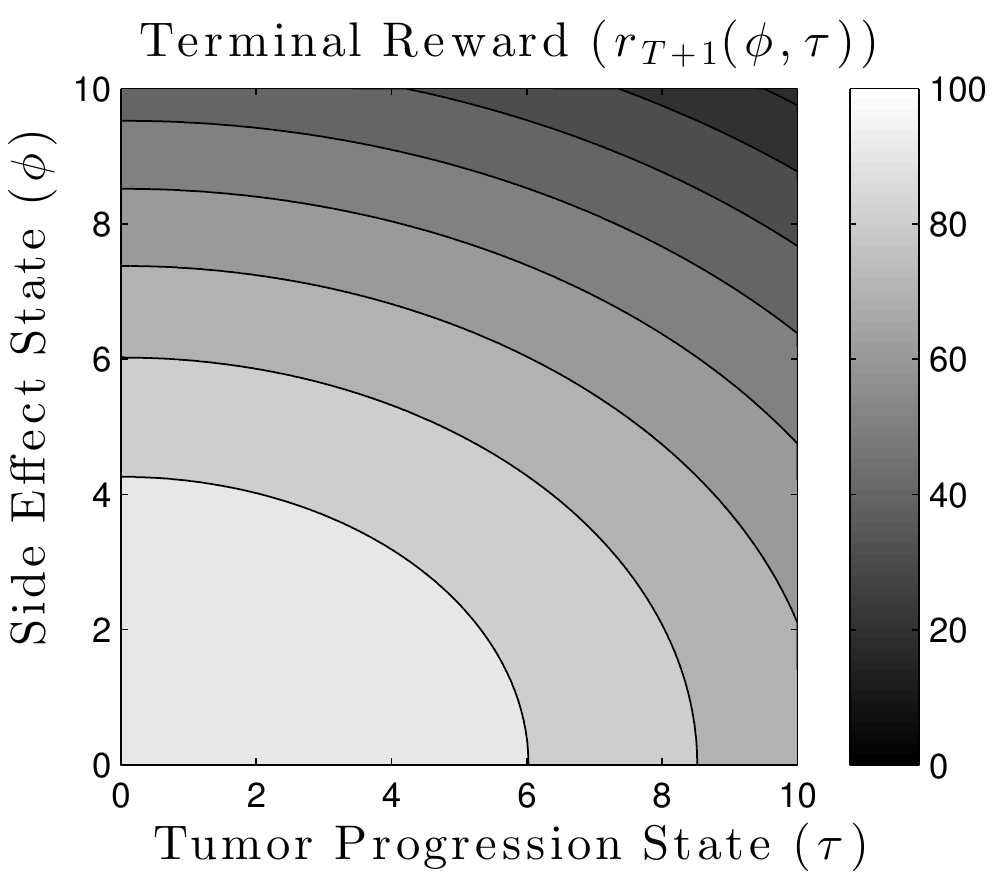} \quad
\includegraphics[height=2in]{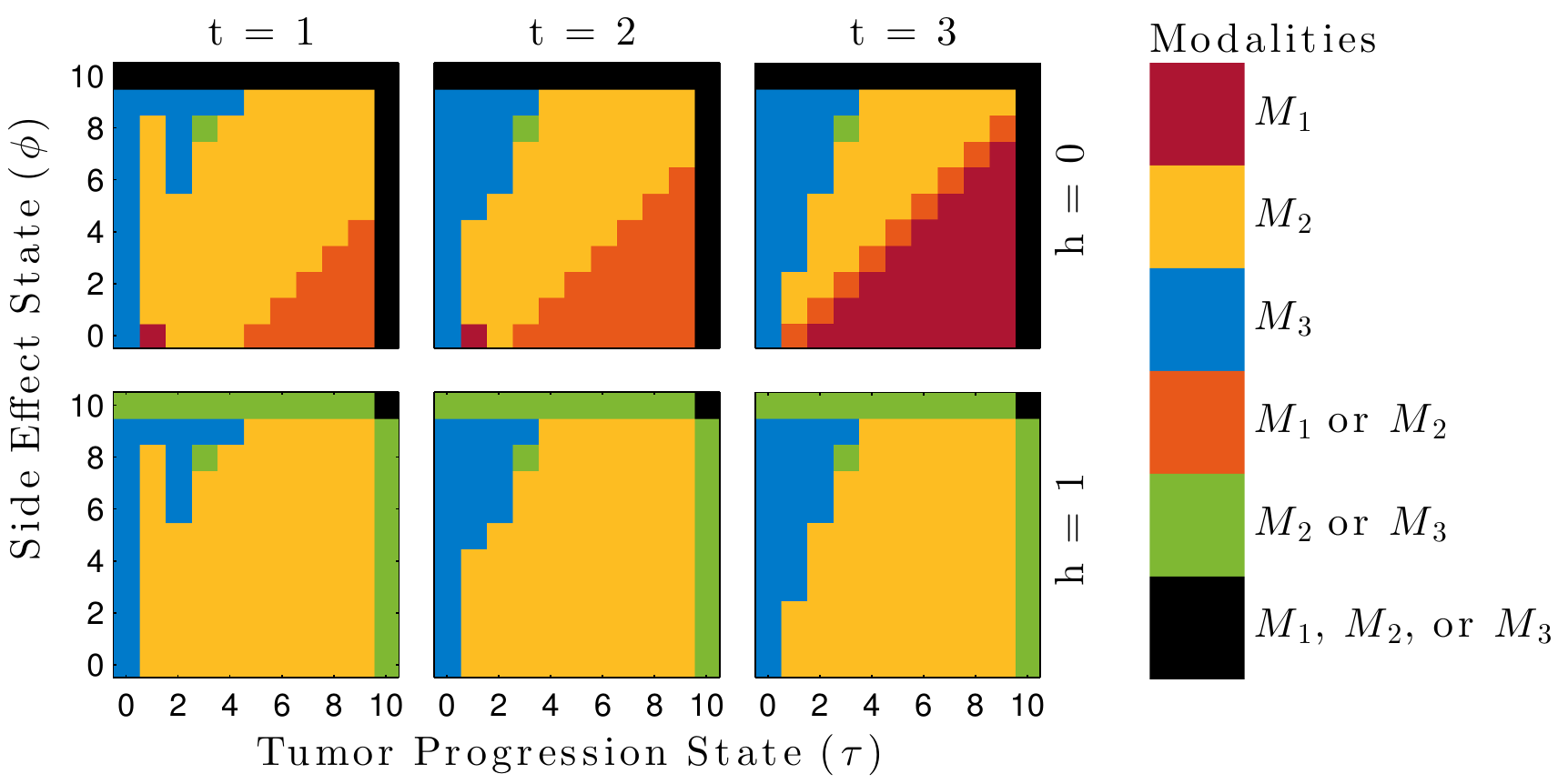}
\caption{(Left) Quadratic terminal reward function $r_{T+1}(\phi,\tau) = \frac{1}{3} f(\phi; 2) + \frac{2}{3} g(\tau; 2)$. (Right) Optimal treatment policy for $T=3$ and $A=\{M_1, M_2, M_3\}$ with transition probabilities given in Table \ref{tab:baseprobability}.}
\label{fig:c33}
\end{figure}

\begin{figure}[H]
\captionsetup{width=0.9\textwidth}
\centering
\includegraphics[height=2in]{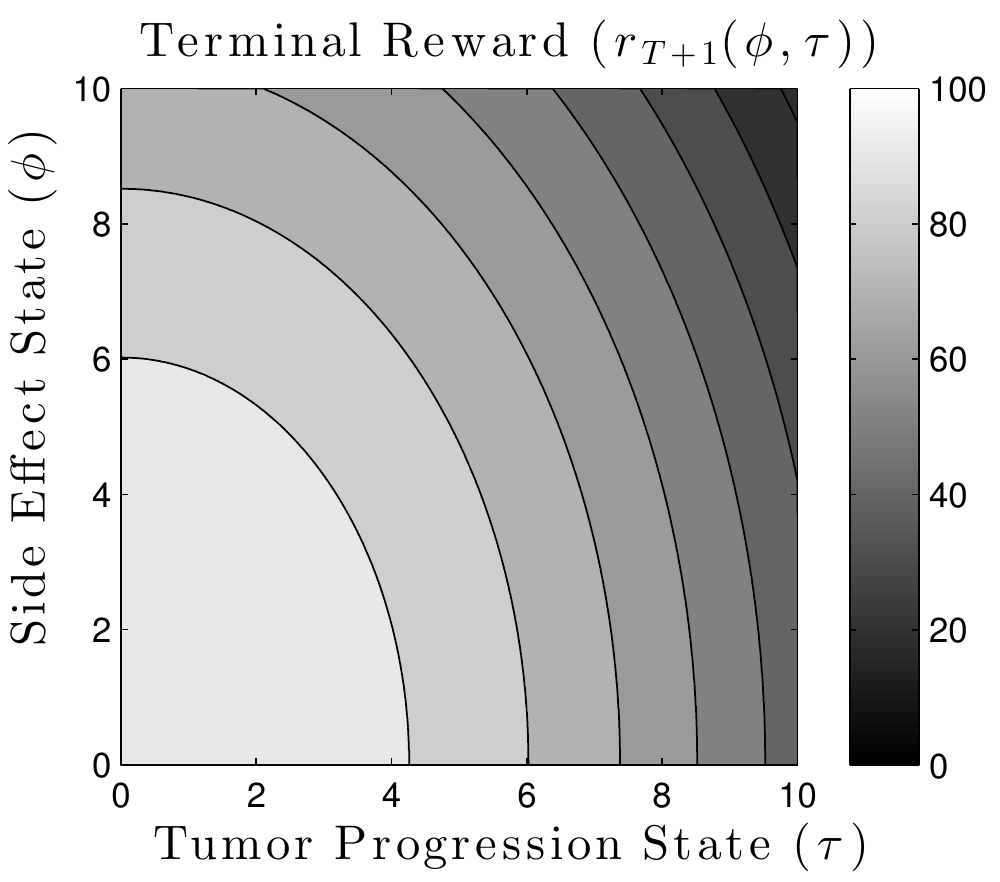} \quad
\includegraphics[height=2in]{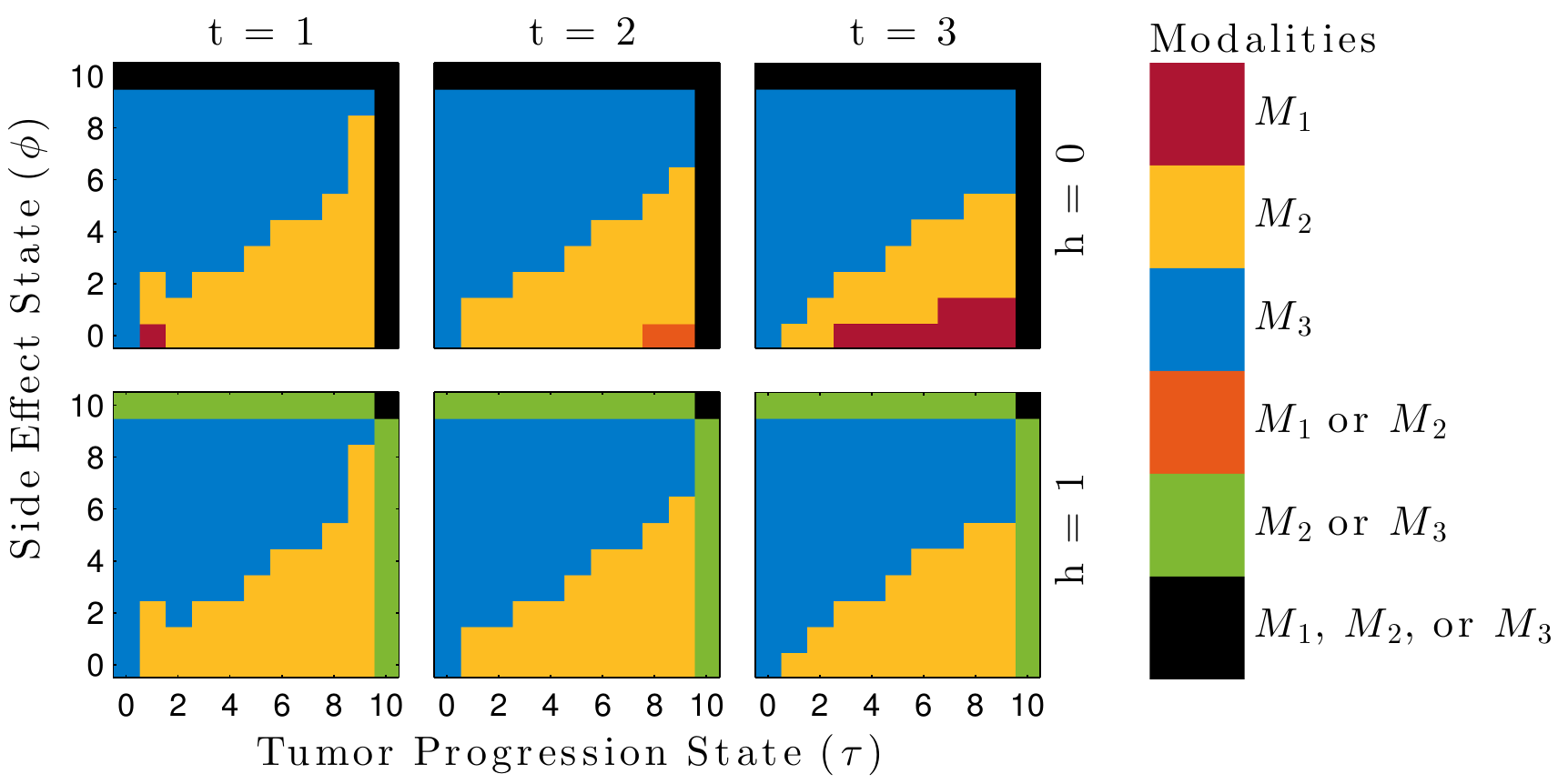}
\caption{(Left) Quadratic terminal reward function $r_{T+1}(\phi,\tau) = \frac{2}{3} f(\phi; 2) + \frac{1}{3} g(\tau; 2)$. (Right) Optimal treatment policy for $T=3$ and $A=\{M_1, M_2, M_3\}$ with transition probabilities given in Table \ref{tab:baseprobability}.}
\label{fig:c67}
\end{figure}

Next we look at the effect of the relative importance of side effect and tumor progression in the concave terminal reward function with $c_\phi = c$ and $c_\tau = (1-c)$, given by
\begin{eqnarray}
r_{T+1}(\phi,\tau) \hspace{-0.8em}&=&\hspace{-0.8em} c f(\phi;2) + (1-c) g(\tau;2) \\
\hspace{-0.8em}&=&\hspace{-0.8em} c \Big(100-\phi^2\Big) + (1-c) \left(100-\tau^2\right),
\end{eqnarray}
where $c \in \{1/3, 1/2 (\mbox{base case}), 2/3\}$. The optimal policies with $c=1/3$ and $2/3$ are shown in Figure \ref{fig:c33} and Figure \ref{fig:c67} respectively.  

As we weight tumor progression more heavily than side effect in the terminal reward function, the optimal policies become more aggressive. Specifically, decreasing $c$ produces policies that recommend treatment over surveillance in more states with an increased frequency of using $M_1$ when $h=0$ and $M_2$ when $h=1$. 

\subsection{Effect of intermediate rewards on optimal policies}
\label{subsec:intermediate}

In this section we explore the effect of adding an intermediate reward function on the optimal policies. We consider two cases, where we collect rewards for (1) reducing side effect during the treatment course by adding $r_\phi$, and (2) reducing tumor progression during the treatment course by adding $r_\tau$.  Using quadratic rewards ($d_\phi=d_\tau=2$), the intermediate reward functions $r_\phi$ and $r_\tau$ are defined as
\begin{eqnarray}
r_\phi(s_t,s_{t+1}) \hspace{-0.8em}& =&\hspace{-0.8em} c_m f(\phi_{t+1}; 2), \quad t=1,2,\dots,T, \\
r_\tau(s_t,s_{t+1}) \hspace{-0.8em}&= &\hspace{-0.8em} c_m g(\tau_{t+1}; 2), \quad\, t=1,2,\dots,T,
\end{eqnarray} 
where we let $c_m=1/4$.

\begin{figure}[H]
\captionsetup{width=0.9\textwidth}
\centering
\includegraphics[height=2in]{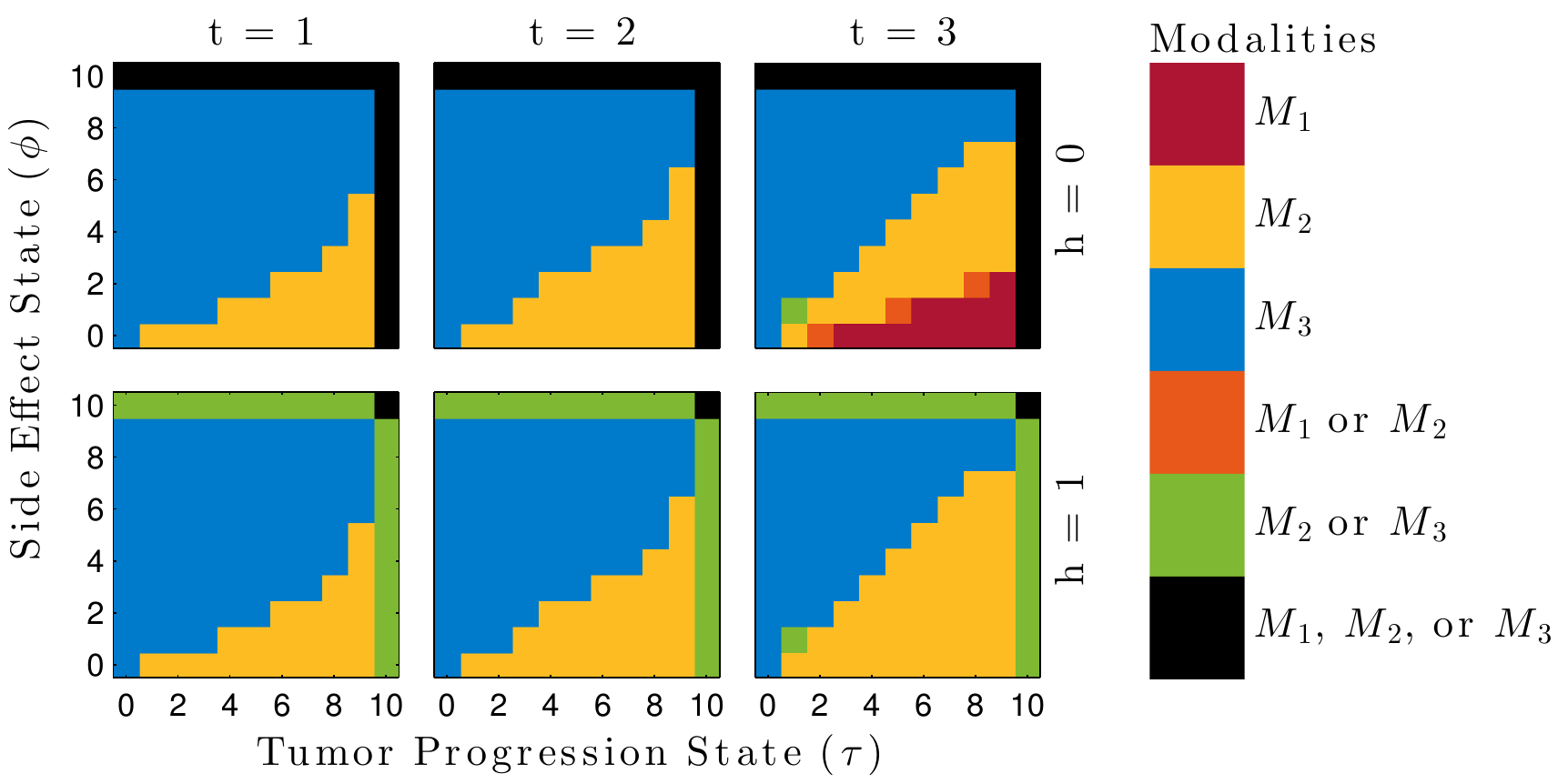}
\caption{Optimal treatment policy for $T=3$ and $A=\{M_1, M_2, M_3\}$ with transition probabilities given in Table \ref{tab:baseprobability}, terminal reward function $r_{T+1}(\phi,\tau) = \frac{1}{2} f(\phi; 2) + \frac{1}{2} g(\tau; 2)$, and intermediate reward function $r_\phi(\phi) = \frac{1}{4} f(\phi; 2)$.}
\label{fig:intermediateSideeffect}
\end{figure}

\begin{figure}[H]
\captionsetup{width=0.9\textwidth}
\centering
\includegraphics[height=2in]{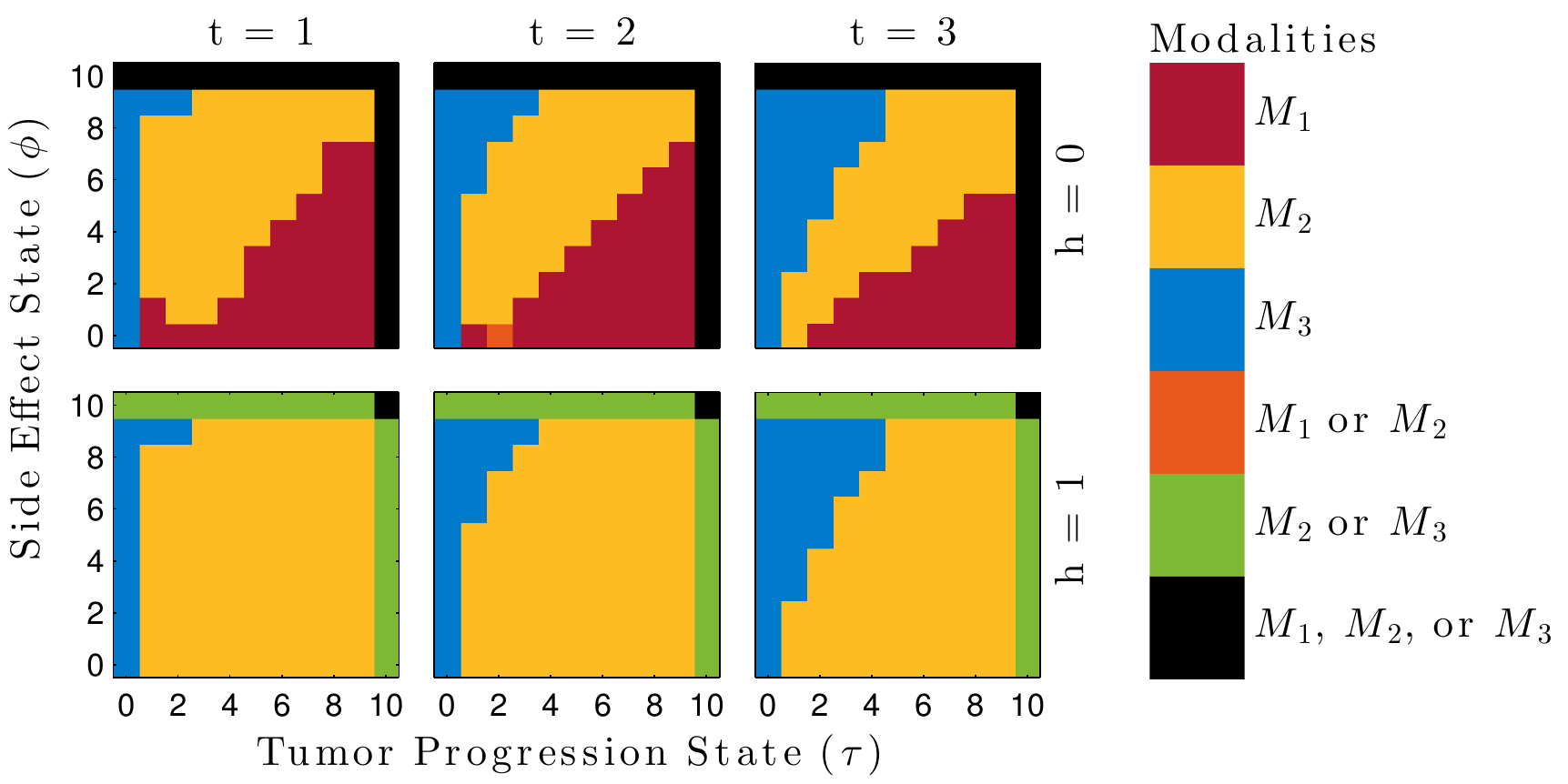}
\caption{Optimal treatment policy for $T=3$ and $A=\{M_1, M_2, M_3\}$ with transition probabilities given in Table \ref{tab:baseprobability}, terminal reward function $r_{T+1}(\phi,\tau) = \frac{1}{2} f(\phi; 2) + \frac{1}{2} g(\tau; 2)$, and intermediate reward function $r_\tau(\tau) = \frac{1}{4} g(\tau; 2)$.}
\label{fig:intermediateTumor}
\end{figure}

The optimal policies using $r_\phi$ and $r_\tau$ are presented in Figure \ref{fig:intermediateSideeffect} and Figure \ref{fig:intermediateTumor} respectively.  

The optimal policy using $r_\phi$ tends to select actions that reduce side effect during the treatment course; therefore,  the optimal policy has a higher proportion of surveillance ($M_3$) and a lower proportion of $M_1$  and $M_2$ with $h=0$, and $M_2$ with $h=1$ in all treatment periods compared to the base case in Figure \ref{fig:basecase}. We note that $M_1$, which can only be used once, is saved for the last treatment period ($t=T$) when $r_\phi$ is used. This can be compared with Figure \ref{fig:c67}, where the importance of reducing side effect is larger in the terminal reward than the importance of reducing tumor progression. When $c_\phi$ is larger than $c_\tau$ without intermediate rewards, the optimal policy still includes $M_1$ for certain states in earlier treatment periods.

On the other hand, the optimal policy using $r_\tau$ tends to select actions that reduce tumor progression during the treatment course.  In this case, the proportion of $M_3$ decreases in all treatment periods and the proportion of $M_1$ and $M_2$ with $h=0$, and $M_2$ with $h=1$ increases. Therefore the preference for treatment over surveillance is higher.  We note the difference in the optimal policy resulting from adding $r_\tau$ to that of Figure \ref{fig:c33}, where $c_\tau$ is increased in the terminal reward. The proportion of $M_1$, which can only be used once, is higher in all treatment periods when $r_\tau$ is used.  This is different from using a larger $c_\tau$ than $c_\phi$ in the terminal reward function, where $M_1$ is used in a larger portion of the states in the last treatment period than in the earlier periods.

\subsection{Effect of transition probabilities on optimal policies}
\label{subsec:varioustransition}

In this section we show how changes in transition probabilities affect the optimal treatment policy. First, we explore a case where the treatment modality $M_1$ is more effective in reducing tumor progression than in the base case, where we increase the probability that the tumor progression will decrease after treatment, that is, 
\begin{equation}
P^{\mathcal{T}}(\tau_t - 1 | \tau_t, M_1) = 0.8 \quad \text{and} \quad P^{\mathcal{T}}(\tau_t | \tau_t, M_1) = 0.2.
\end{equation}

\begin{table}[H]
\captionsetup{width=0.9\textwidth}
\centering
\begin{tabularx}{0.9\textwidth}{c *{6}{Y}}
\toprule
\multirow{2}{*}{\bf Modality} 
& \multicolumn{3}{c}{\bf Side Effect in period $t+1$}  
& \multicolumn{3}{c}{\bf Tumor Progression in period $t+1$}\\
\multirow{2}{*}{$(a_t)$} 
& \multicolumn{3}{c}{$(\phi_{t+1})$} 
& \multicolumn{3}{c}{$(\tau_{t+1})$} \\ 
\cmidrule(lr){2-4} \cmidrule(l){5-7}
& $\phi_t-1$ & $\phi_t$ & $\phi_t+1$ & $\tau_t-1$ & $\tau_t$ & $\tau_t+1$ \\
\midrule
$M_1$ &    0 & 0.4 & 0.6 & \bf{ 0.8} & \bf{0.2} &    0 \\
$M_2$ &    0 & 0.6 & 0.4 & 0.6 & 0.4 &    0 \\
$M_3$ & 0.6 & 0.4 &    0 &    0 & 0.3 & 0.7 \\
\bottomrule
\end{tabularx}
\caption{State transition probabilities $P^{\Phi}(\phi_{t+1}|\phi_t,a_t)$ and $P^{\mathcal{T}}(\tau_{t+1}|\tau_t,a_t)$ with the assumption that $M_1$ is more effective than the base case in reducing tumor progression.}
\label{tab:probability2}
\end{table}

\begin{figure}[H]
\captionsetup{width=0.9\textwidth}
\centering
\includegraphics[height=2in]{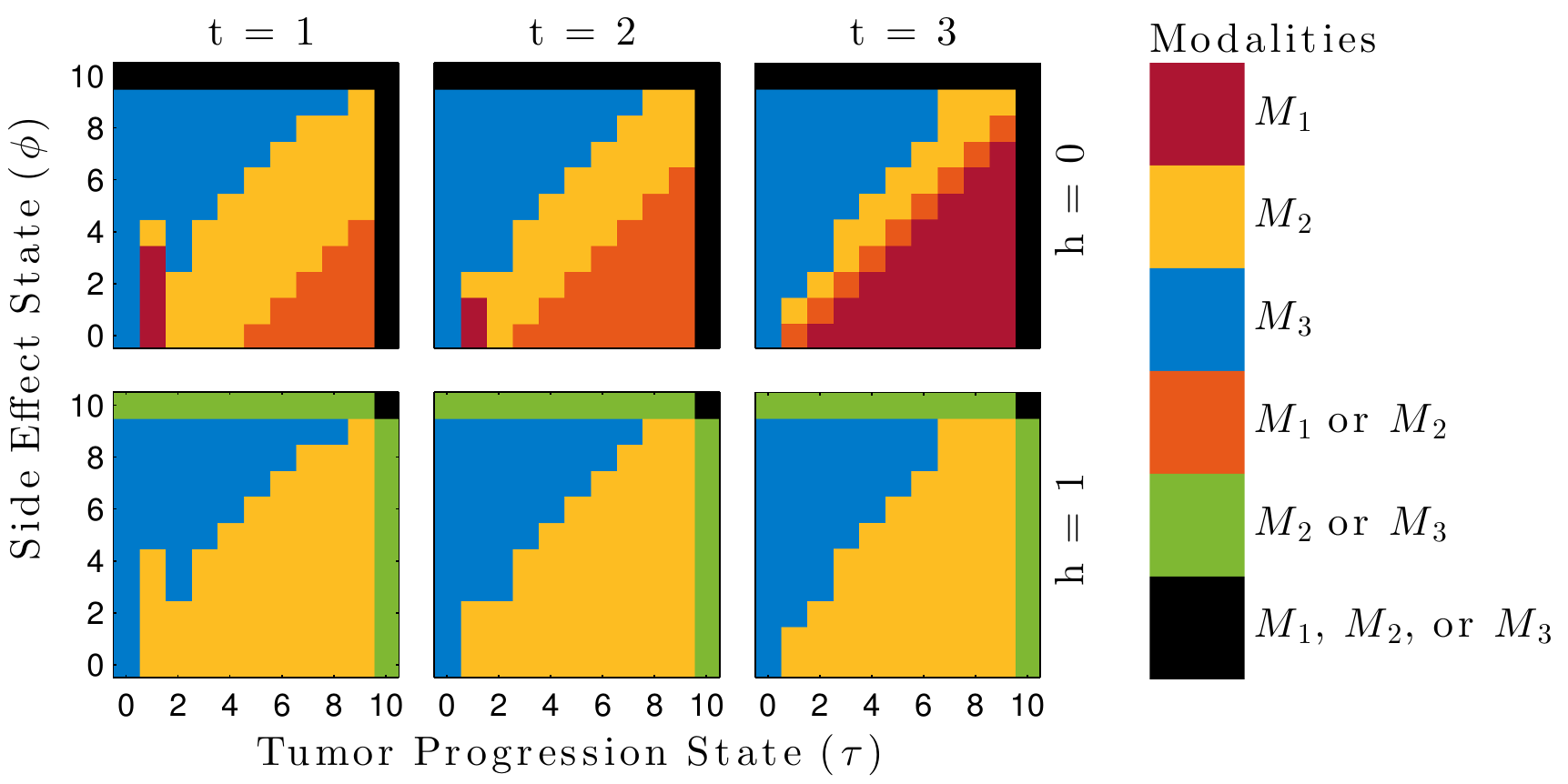}
\caption{Optimal treatment policy for $T=3$ and $A=\{M_1, M_2, M_3\}$ with transition probabilities given in Table \ref{tab:probability2}, modified to make $M_1$ more effective than the base case in decreasing tumor progression, and quadratic terminal reward function $r_{T+1}(\phi,\tau) = \frac{1}{2} f(\phi,2) + \frac{1}{2} g(\tau,2)$.}
\label{fig:prob2}
\end{figure}

The optimal policy computed using the transition probabilities in Table \ref{tab:probability2} is shown in Figure \ref{fig:prob2}. It shows a more aggressive treatment policy that uses $M_1$ in a higher proportion of the states compared to the base case.  We note an interesting effect of the boundary condition at $\tau = 0$, specifically, $M_1$ is chosen for several states where $\tau = 1$ and $\phi < 4$ (without severe side effect) at $t=1$ due to the fact that tumor remission in earlier treatment periods ($t=1, 2$) is now more likely if $M_1$ is chosen.

Next we explore a case where the treatment modality $M_2$ is less risky than that in the base case in terms of increasing side effect, where we decrease the probability that the side effect will increase after treatment, that is,
\begin{equation}
P^{\Phi}(\phi_t | \phi_t, M_2) = 0.7 \quad \text{and} \quad P^{\Phi}(\phi_t+1 | \phi_t, M_2) = 0.3.
\end{equation}

\begin{figure}[H]
\captionsetup{width=0.9\textwidth}
\centering
\includegraphics[height=2in]{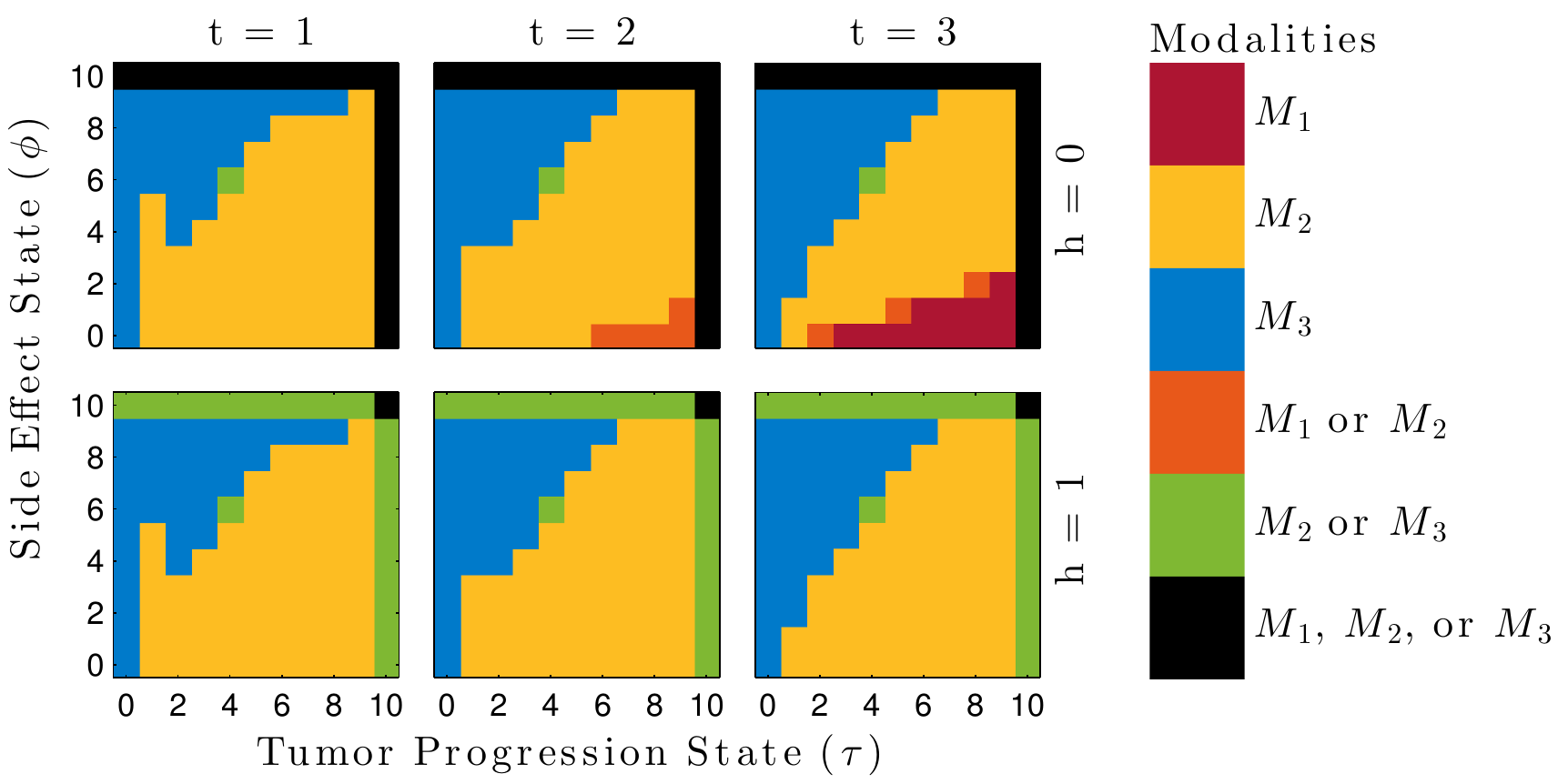}
\caption{Optimal treatment policy for $T=3$ and $A=\{M_1, M_2, M_3\}$ with transition probabilities given in Table \ref{tab:probability3}, modified to make $M_2$ less risky than the base case in increasing side effect, and quadratic terminal reward function $r_{T+1}(\phi,\tau) = \frac{1}{2} f(\phi,2) + \frac{1}{2} g(\tau,2)$.}
\label{fig:prob3}
\end{figure}

The state transition probabilities in this case are shown in Table \ref{tab:probability3} and the resulting optimal policy is shown in Figure \ref{fig:prob3}. We see that this change makes $M_2$ more favorable and produces an optimal policy that uses $M_2$ in a higher proportion of the states in all treatment periods. 

\begin{table}[H]
\captionsetup{width=0.9\textwidth}
\centering
\begin{tabularx}{0.9\textwidth}{c *{6}{Y}}
\toprule
\multirow{2}{*}{\bf Modality} 
& \multicolumn{3}{c}{\bf Side Effect in period $t+1$}  
& \multicolumn{3}{c}{\bf Tumor Progression in period $t+1$}\\
\multirow{2}{*}{$(a_t)$} 
& \multicolumn{3}{c}{$(\phi_{t+1})$} 
& \multicolumn{3}{c}{$(\tau_{t+1})$} \\ 
\cmidrule(lr){2-4} \cmidrule(l){5-7}
& $\phi_t-1$ & $\phi_t$ & $\phi_t+1$ & $\tau_t-1$ & $\tau_t$ & $\tau_t+1$ \\
\midrule
$M_1$ &    0 & 0.4 & 0.6 & 0.7 & 0.3 &    0 \\
$M_2$ &    0 & \bf{0.7} & \bf{0.3} & 0.6 & 0.4 &    0 \\
$M_3$ & 0.6 & 0.4 &    0 &    0 & 0.3 & 0.7 \\
\bottomrule
\end{tabularx}
\caption{State transition probabilities $P^{\Phi}(\phi_{t+1}|\phi_t,a_t)$ and $P^{\mathcal{T}}(\tau_{t+1}|\tau_t,a_t)$ with the assumption that $M_2$ is less risky than the base case in terms of increasing side effect.}
\label{tab:probability3}
\end{table}

\begin{figure}[H]
\captionsetup{width=0.9\textwidth}
\centering
\includegraphics[height=2in]{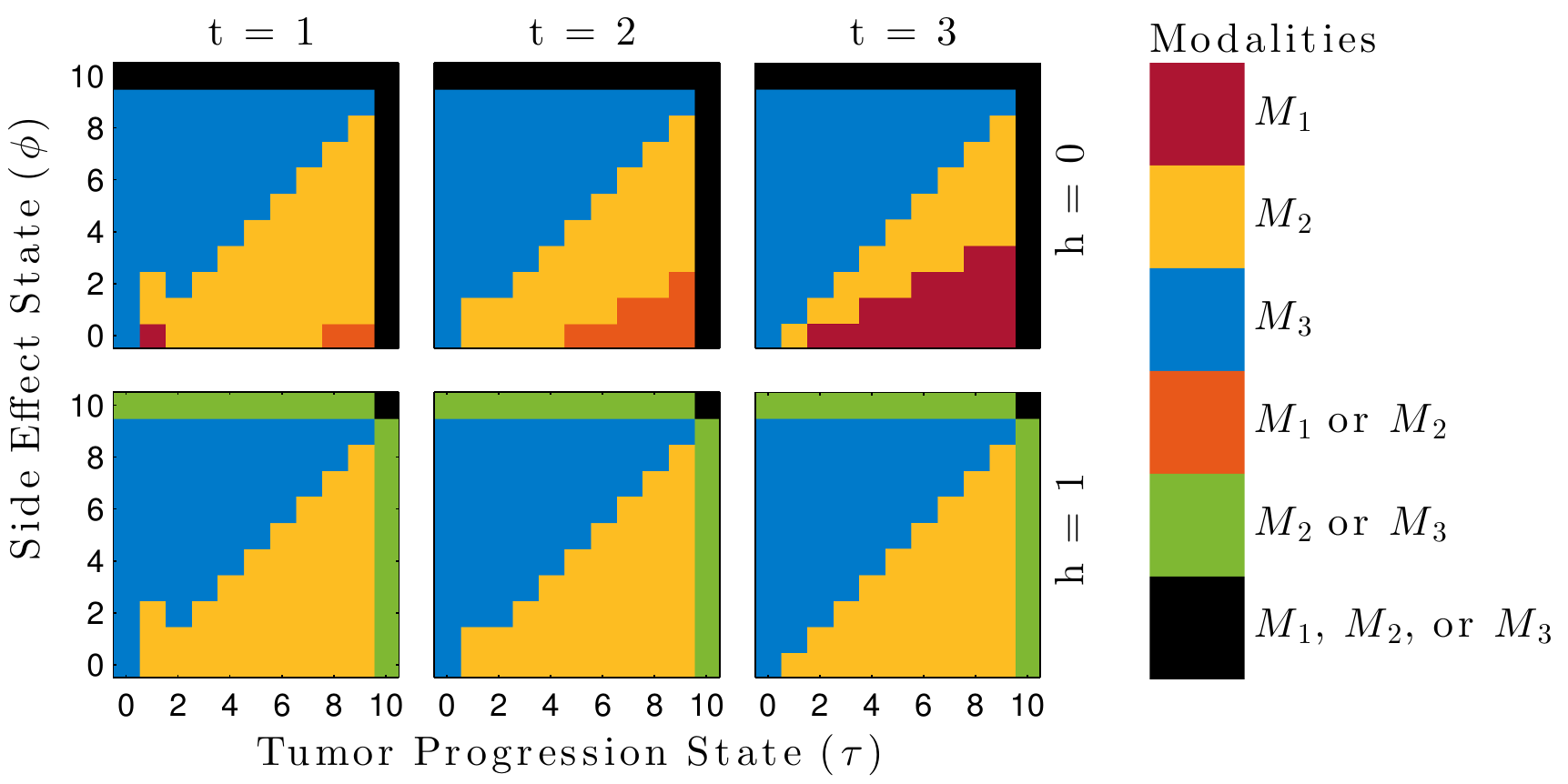}
\caption{Optimal treatment policy for $T=3$ and $A=\{M_1, M_2, M_3\}$ with transition probabilities given in Table \ref{tab:probability4}, modified to make $M_3$ less risky than the base case in increasing tumor progression, and quadratic terminal reward function $r_{T+1}(\phi,\tau) = \frac{1}{2} f(\phi,2) + \frac{1}{2} g(\tau,2)$.}
\label{fig:prob4}
\end{figure}

\newpage
Finally, we investigate a case where the surveillance modality $M_3$ is less likely to increase tumor progression than the base case, where we change probabilities in the fourth row of Table \ref{tab:baseprobability} to
\begin{equation}
P^{\mathcal{T}}(\tau_t  | \tau_t, M_3) = 0.7 \quad \text{and} \quad P^{\mathcal{T}}(\tau_t+1 | \tau_t, M_3) = 0.3.
\end{equation}
The state transition probabilities in this case are shown in Table \ref{tab:probability4} and the resulting optimal policy is shown in Figure \ref{fig:prob4}. We see an increase in the amount of surveillance suggested by the policy in all treatment periods. We also note that $M_3$ is optimal in more states with higher side effects compared to the base case. This policy may be appropriate in the case of a slowly growing tumor, where the probability of increasing tumor progression by surveillance is relatively low for a given treatment time period.

\begin{table}[H]
\captionsetup{width=0.9\textwidth}
\centering
\begin{tabularx}{0.9\textwidth}{c *{6}{Y}}
\toprule
\multirow{2}{*}{\bf Modality} 
& \multicolumn{3}{c}{\bf Side Effect in period $t+1$}  
& \multicolumn{3}{c}{\bf Tumor Progression in period $t+1$}\\
\multirow{2}{*}{$(a_t)$} 
& \multicolumn{3}{c}{$(\phi_{t+1})$} 
& \multicolumn{3}{c}{$(\tau_{t+1})$} \\ 
\cmidrule(lr){2-4} \cmidrule(l){5-7}
& $\phi_t-1$ & $\phi_t$ & $\phi_t+1$ & $\tau_t-1$ & $\tau_t$ & $\tau_t+1$ \\
\midrule
$M_1$ &    0 & 0.4 & 0.6 & 0.7 & 0.3 &    0 \\
$M_2$ &    0 & 0.6 & 0.4 & 0.6 & 0.4 &    0 \\
$M_3$ & 0.6 & 0.4 &    0 &    0 & \bf{0.7} & \bf{0.3} \\
\bottomrule
\end{tabularx}
\caption{State transition probabilities $P^{\Phi}(\phi_{t+1}|\phi_t,a_t)$ and $P^{\mathcal{T}}(\tau_{t+1}|\tau_t,a_t)$ with the assumption that $M_2$ is less risky than the base case in terms of increasing side effect.}
\label{tab:probability4}
\end{table}

\subsection{Multiple Type 2 modalities}
\label{subsec:multipleM2}

In this section we demonstrate our model with multiple modalities of a particular treatment type. Specifically, we consider a case where there are two modalities of Type 2, i.e., $A=\{ M_1, M_{2a}, M_{2b}, M_3 \}$. We order the index of the treatment modalities based on effectiveness, that is,
\begin{align}
P^\Phi(\phi_{t+1}|\phi_t,M_1) \geq P^\Phi(\phi_{t+1}|\phi_t,M_{2a}) &\geq P^\Phi(\phi_{t+1}|\phi_t,M_{2b}) \geq P^\Phi(\phi_{t+1}|\phi_t,M_{3}) \text{ for } \phi_{t+1} > \phi_t,\quad  \\
P^\mathcal{T}(\tau_{t+1}|\tau_t,M_1) \geq P^\mathcal{T}(\tau_{t+1}|\tau_t,M_{2a}) &\geq P^\mathcal{T}(\tau_{t+1}|\tau_t,M_{2b}) \geq P^\mathcal{T}(\tau_{t+1}|\tau_t,M_{3}) \text{ for } \tau_{t+1} < \tau_t.
\end{align}

\begin{figure}[H]
\captionsetup{width=0.9\textwidth}
\centering
\includegraphics[height=2in]{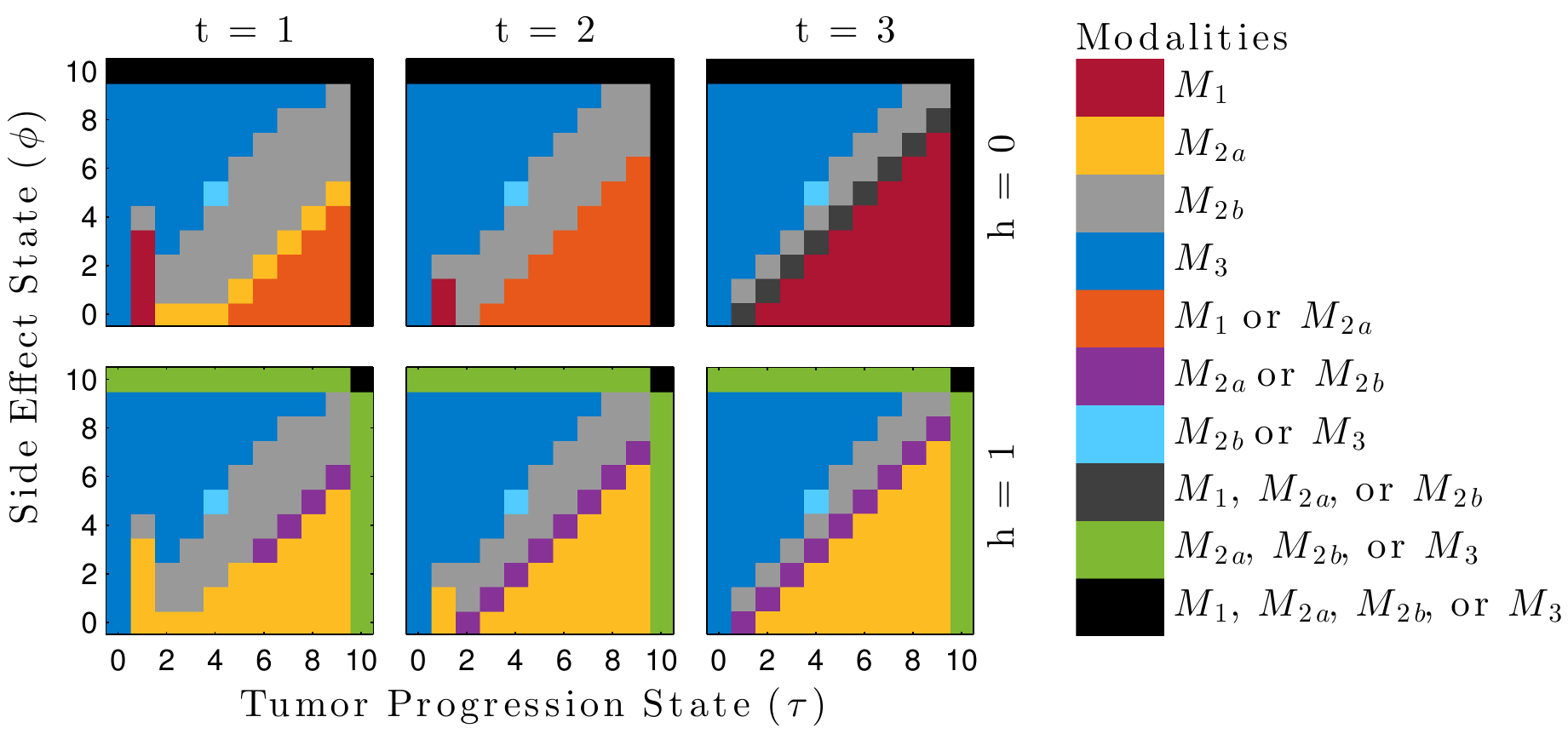}
\caption{Optimal treatment policy for $T=3$ and $A=\{M_1, M_{2a}, M_{2b}, M_3\}$ with probabilities given in Table \ref{tab:4probability} and quadratic terminal reward function $r_{T+1}(\phi,\tau) = \frac{1}{2} f(\phi;2) + \frac{1}{2} g(\tau; 2)$.}
\label{fig:4actions}
\end{figure}

The transition probabilities used in this case are presented in Table \ref{tab:4probability}, and the resulting optimal policy is shown in Figure \ref{fig:4actions}. The trend is qualitatively similar to the previous cases with three modalities.  We note that there are only ties between consecutive treatment modalities, for example, $M_1$ never ties with $M_{2b}$ or $M_3$, because of the ordering in the action space. 

\begin{table}[H]
\captionsetup{width=0.9\textwidth}
\centering
\begin{tabularx}{0.9\textwidth}{c *{6}{Y}}
\toprule
\multirow{2}{*}{\bf Modality} 
& \multicolumn{3}{c}{\bf Side Effect in period $t+1$}  
& \multicolumn{3}{c}{\bf Tumor Progression in period $t+1$}\\
\multirow{2}{*}{$(a_t)$} 
& \multicolumn{3}{c}{$(\phi_{t+1})$} 
& \multicolumn{3}{c}{$(\tau_{t+1})$} \\ 
\cmidrule(lr){2-4} \cmidrule(l){5-7}
& $\phi_t-1$ & $\phi_t$ & $\phi_t+1$ & $\tau_t-1$ & $\tau_t$ & $\tau_t+1$ \\
\midrule
$M_1$     &   0 & 0.4 & 0.6 & 0.7 & 0.3 &    0 \\
$M_{2a}$ &   0 & 0.5 & 0.5 & 0.6 & 0.4 &    0 \\
$M_{2b}$ &   0 & 0.6 & 0.4 & 0.5 & 0.5 &    0 \\
$M_3$     & 0.6 & 0.4 &    0 &    0 & 0.3 & 0.7 \\
\bottomrule
\end{tabularx}
\caption{State transition probabilities $P^{\Phi}(\phi_{t+1}|\phi_t,a)$ and $P^{\mathcal{T}}(\tau_{t+1}|\tau_t,a)$ used for a case with $T=3$ and \\$A = \{M_1, M_{2a}, M_{3a}, M_4\}$}
\label{tab:4probability}
\end{table}

\section{Conclusions and future research}
\label{sec:conclusion} 

With diverse patient characteristics, numerous treatment modalities available in modern medicine, and various possible outcomes, making treatment decisions tailored to each individual patient is extremely complex.  It may no longer be practical to make optimal decisions based solely on individual clinician's experiences and empirical intuition. We proposed a novel mathematical framework to model optimal treatment policies for cancer therapy using a finite-horizon Markov decision process (MDP). Numerical simulations using simplified patient states and clinically intuitive reward functions have shown the potential application of our model to aid in treatment decision-making. Using state transition probabilities obtained from treatment-outcome clinical data, our model can assist clinicians in making optimal decisions for the patient's current state.  Our model has the potential to be personalized to individual patients via custom utility functions based on the patient's own preferences, with the addition of state variables such as age, specific side effects of interests, and tumor types. It can also serve as a tool to explain an expected treatment course to patients.

As we take more realistic and detailed attributes into consideration in the model, the problem becomes computationally intractable due to increased state and outcome space.  We leave these high dimensional problems, which require an approximate dynamic programming approach, for future work.  Other potential formulations for future investigation include infinite-horizon MDPs leading to stationary optimal solutions. An infinite horizon formulation is particularly useful when the patient's expected treatment course is unpredictable.

\section*{Acknowledgement}
This study is supported in part by National Science Foundation grant CMMI 1560476.

\bibliographystyle{unsrt}
\bibliography{references}

\begin{thebibliography}{10}

\bibitem{siegel2016cancerfacts}
Rebecca Siegel, Kimberly Miller, and Ahmedin Jemal.
\newblock Cancer facts \& figures, 2016.
\newblock {\em Atlanta: American Cancer Society}, pages 1--66, 2016.

\bibitem{torre2015global}
Lindsey Torre, Rebecca Siegel, and Ahmedin Jemal.
\newblock Global cancer facts \& figures, 3rd edition.
\newblock {\em Atlanta: American Cancer Society}, pages 1--57, 2015.

\bibitem{siegel2016cancerstats}
Rebecca~L Siegel, Kimberly~D Miller, and Ahmedin Jemal.
\newblock Cancer statistics, 2016.
\newblock {\em CA: A Cancer Journal for Clinicians}, 66(1):7--30, 2016.

\bibitem{kemp2016combo}
Jessica~A Kemp, Min~Suk Shim, Chan~Yeong Heo, and Young~Jik Kwon.
\newblock ``{C}ombo'' nanomedicine: Co-delivery of multi-modal therapeutics for
  efficient, targeted, and safe cancer therapy.
\newblock {\em Advanced Drug Delivery Reviews}, 98:3--18, 2016.

\bibitem{franco2016combined}
Pierfrancesco Franco, Alba Fiorentino, Francesco Dionisi, Michele Fiore, Silvia
  Chiesa, Stefano Vagge, Francesco Cellini, Luciana Caravatta, Mario Tombolini,
  Fiorenza De~Rose, et~al.
\newblock Combined modality therapy for thoracic and head and neck cancers: A
  review of updated literature based on a consensus meeting.
\newblock {\em Tumori}, 102(5):459, 2016.

\bibitem{sathyanarayanan2015cancer}
Vishwanath Sathyanarayanan and Sattva~S Neelapu.
\newblock Cancer immunotherapy: Strategies for personalization and
  combinatorial approaches.
\newblock {\em Molecular Oncology}, 9(10):2043--2053, 2015.

\bibitem{taunk2016external}
Neil~K Taunk, Fabio~Y Moraes, Freddy~E Escorcia, Lucas~Castro Mendez, Kathryn
  Beal, and Gustavo~N Marta.
\newblock External beam re-irradiation, combination chemoradiotherapy, and
  particle therapy for the treatment of recurrent glioblastoma.
\newblock {\em Expert Review of Anticancer Therapy}, 16(3):347--358, 2016.

\bibitem{matzenauer2016treatment}
Marcel Matzenauer, David Vrana, Bohuslav Melichar, et~al.
\newblock Treatment of brain metastases.
\newblock {\em Biomedical Papers}, 160(4):484--490, 2016.

\bibitem{artac2016update}
Mehmet Artac, Levent Korkmaz, Bassel El-Rayes, and Philip~A Philip.
\newblock An update on the multimodality of localized rectal cancer.
\newblock {\em Critical Reviews in Oncology/Hematology}, 108:23--32, 2016.

\bibitem{shapiro2015neoadjuvant}
Joel Shapiro, J~Jan~B Van~Lanschot, Maarten~CCM Hulshof, Pieter van Hagen,
  Mark~I van Berge~Henegouwen, Bas~PL Wijnhoven, Hanneke~WM van Laarhoven,
  Grard~AP Nieuwenhuijzen, Geke~AP Hospers, Johannes~J Bonenkamp, et~al.
\newblock Neoadjuvant chemoradiotherapy plus surgery versus surgery alone for
  oesophageal or junctional cancer (cross): Long-term results of a randomised
  controlled trial.
\newblock {\em The Lancet Oncology}, 16(9):1090--1098, 2015.

\bibitem{schaefer2005modeling}
Andrew~J Schaefer, Matthew~D Bailey, Steven~M Shechter, and Mark~S Roberts.
\newblock Modeling medical treatment using {M}arkov decision processes.
\newblock In {\em Operations Research and Health Care}, pages 593--612.
  Springer, 2005.

\bibitem{bennett2013artificial}
Casey~C Bennett and Kris Hauser.
\newblock Artificial intelligence framework for simulating clinical
  decision-making: A {M}arkov decision process approach.
\newblock {\em Artificial Intelligence in Medicine}, 57(1):9--19, 2013.

\bibitem{beil2001analysis}
Damian~R Beil and Lawrence~M Wein.
\newblock Analysis and comparison of multimodal cancer treatments.
\newblock {\em Mathematical Medicine and Biology}, 18(4):343--376, 2001.

\bibitem{hathout2016modeling}
Leith Hathout, Benjamin Ellingson, and Whitney Pope.
\newblock Modeling the efficacy of the extent of surgical resection in the
  setting of radiation therapy for glioblastoma.
\newblock {\em Cancer Science}, 107(8):1110--1116, 2016.

\bibitem{roy1999coastal}
Nicholas Roy and Sebastian Thrun.
\newblock Coastal navigation with mobile robots.
\newblock In {\em NIPS}, volume~13, pages 1043--1049, 1999.

\bibitem{kinjo2015evaluation}
Ken Kinjo, Eiji Uchibe, and Kenji Doya.
\newblock Evaluation of linearly solvable {M}arkov decision process with
  dynamic model learning in a mobile robot navigation task.
\newblock {\em Value and Reward Based Learning in Neurobots}, page~75, 2015.

\bibitem{briggs1998introduction}
Andrew Briggs and Mark Sculpher.
\newblock An introduction to {M}arkov modelling for economic evaluation.
\newblock {\em Pharmacoeconomics}, 13(4):397--409, 1998.

\bibitem{rust1996numerical}
John Rust.
\newblock Numerical dynamic programming in economics.
\newblock {\em Handbook of Computational Economics}, 1:619--729, 1996.

\bibitem{magni2000deciding}
Paolo Magni, Silvana Quaglini, Monia Marchetti, and Giovanni Barosi.
\newblock Deciding when to intervene: A {M}arkov decision process approach.
\newblock {\em International Journal of Medical Informatics}, 60(3):237--253,
  2000.

\bibitem{hauskrecht2000planning}
Milos Hauskrecht and Hamish Fraser.
\newblock Planning treatment of ischemic heart disease with partially
  observable {M}arkov decision processes.
\newblock {\em Artificial Intelligence in Medicine}, 18(3):221--244, 2000.

\bibitem{kim2009markov}
Minsun Kim, Archis Ghate, and MH~Phillips.
\newblock A {M}arkov decision process approach to temporal modulation of dose
  fractions in radiation therapy planning.
\newblock {\em Physics in Medicine and Biology}, 54(14):4455, 2009.

\bibitem{ahn1996involving}
Jae-Hyeon Ahn and John~C Hornberger.
\newblock Involving patients in the cadaveric kidney transplant allocation
  process: A decision-theoretic perspective.
\newblock {\em Management Science}, 42(5):629--641, 1996.

\bibitem{alagoz2004optimal}
Oguzhan Alagoz, Lisa~M Maillart, Andrew~J Schaefer, and Mark~S Roberts.
\newblock The optimal timing of living-donor liver transplantation.
\newblock {\em Management Science}, 50(10):1420--1430, 2004.

\bibitem{nccncns}
{\em {NCCN} clinical practice guidelines in oncology - Central Nervous System
  Cancers}.
\newblock National Comprehensive Cancer Network, 2015.

\bibitem{weinstein2009qalys}
Milton~C Weinstein, George Torrance, and Alistair McGuire.
\newblock Qalys: The basics.
\newblock {\em Value in Health}, 12(s1):S5--S9, 2009.

\bibitem{murray1994quantifying}
Christopher~J Murray.
\newblock Quantifying the burden of disease: The technical basis for
  disability-adjusted life years.
\newblock {\em Bulletin of the World Health Organization}, 72(3):429, 1994.

\bibitem{mehrez1989quality}
Abraham Mehrez and Amiram Gafni.
\newblock Quality-adjusted life years, utility theory, and healthy-years
  equivalents.
\newblock {\em Medical Decision Making}, 9(2):142--149, 1989.

\bibitem{Saokaew2016}
S~Saokaew, A~Rayanakorn, DB~Wu, and N~Chaiyakunapruk.
\newblock Cost effectiveness of peumococcal vaccination in children in
  low-and-middle income countries: a systematic review.
\newblock {\em Pharmacoeconomics}, 34:1211--1225, 2016.

\bibitem{puterman2014markov}
Martin~L Puterman.
\newblock {\em {M}arkov decision processes: Discrete stochastic dynamic
  programming}.
\newblock John Wiley \& Sons, 2014.

\bibitem{powell2011}
Warren~B Powell.
\newblock {\em Approximate dynamic programming: solving the curses of
  dimensionality}.
\newblock Wiley Series in Probability and Statistics. Wiley, 2nd edition, 2011.

\bibitem{attema2016elicitation}
Arthur~E Attema, Werner~BF Brouwer, Olivier l'Haridon, and Jose~Luis Pinto.
\newblock An elicitation of utility for quality of life under prospect theory.
\newblock {\em Journal of Health Economics}, 48:121--134, 2016.

\bibitem{miravitlles2015clinical}
Marc Miravitlles, Alicia Huerta, Manuel Valle, Patricia Garc{\'\i}a-Sidro,
  Carles Forn{\'e}, Carlos Crespo, and Jos{\'e}~Luis L{\'o}pez-Campos.
\newblock Clinical variables impacting on the estimation of utilities in
  chronic obstructive pulmonary disease.
\newblock {\em International Journal of Chronic Obstructive Pulmonary Disease},
  10:367, 2015.

\bibitem{carradice2011modelling}
D~Carradice, FAK Mazari, N~Samuel, V~Allgar, J~Hatfield, and IC~Chetter.
\newblock Modelling the effect of venous disease on quality of life.
\newblock {\em British Journal of Surgery}, 98(8):1089--1098, 2011.

\bibitem{currie2006multivariate}
Craig~J Currie, Christopher~Ll Morgan, Chris~D Poole, Peter Sharplin, Morten
  Lammert, and Phil McEwan.
\newblock Multivariate models of health-related utility and the fear of
  hypoglycaemia in people with diabetes.
\newblock {\em Current Medical Research and Opinion}, 22(8):1523--1534, 2006.

\end{thebibliography}

\end{document}